# Developing Talent from a Supply–Demand Perspective: An Optimization Model for Managers

**Hadi Moheb-Alizadeh** [1] **and Robert B. Handfield** [2,*]

1. Graduate Program in Operations Research, North Carolina State University, Raleigh, NC 27695, USA; hmoheba@ncsu.edu
2. Department of Business Management, College of Management, North Carolina State University, 2806-A Hillsborough St., Upper Level, Campus Box 7229, Raleigh, NC 27695-7229, USA
* Correspondence: rbhandfi@ncsu.edu; Tel.: +1-(919)-515-4674; Fax: +1-(919)-515-2120

**Abstract:** While executives emphasize that human resources (HR) are a firm's biggest asset, the level of research attention devoted to planning talent pipelines for complex global organizational environments does not reflect this emphasis. Numerous challenges exist in establishing human resource management strategies aligned with strategic operations planning and growth strategies. We generalize the problem of managing talent from a supply–demand standpoint through a resource acquisition lens, to an industrial business case where an organization recruits for multiple roles given a limited pool of potential candidates acquired through a limited number of recruiting channels. In this context, we develop an innovative analytical model in a stochastic environment to assist managers with talent planning in their organizations. We apply supply chain concepts to the problem, whereby individuals with specific competencies are treated as unique products. We first develop a multi-period mixed integer nonlinear programming model and then exploit chance-constrained programming to a linearized instance of the model to handle stochastic parameters, which follow any arbitrary distribution functions. Next, we use an empirical study to validate the model with a large global manufacturing company, and demonstrate how the proposed model can effectively manage talents in a practical context. A stochastic analysis on the implemented case study reveals that a reasonable improvement is derived from incorporating randomness into the problem.

## 1. Introduction

Since creation of the phrase *War for Talent* by a group of McKinsey consultants in 1997 [1], the topic of talent management (TM) has received considerable attention from practitioners and academics. Since that time, issue of talent management represents a concern for an increasing number of organizations engaged in development of their human resources. In recent years, the capability of an organization to manage and plan the skills, knowledge, and competencies of its human resources (i.e., having the right people, with the right skills, in the right place and at the right time) has been viewed as a critical factor for organizational success. To achieve success, investing in human resources is a recognized driver for succeeding in reaching growth, profitability and client satisfaction objectives. The implication is that organizations must effectively and efficiently plan and manage individual workers to reach their highest level of potential, especially for those that have a high impact (such as supply chain talent).

In the context of human resource management (HRM), talent management can be defined as a holistic approach to human resource planning aimed at strengthening organizational capability and driving business priorities using a range of HR interventions [2,3]. These interventions include a focus on performance enhancement, career development and succession planning [4]. Research indicates that talent management must be aligned with corporate strategies and culture to be effective [5]. The concept of TM has evolved into common management practices. While originally focused on recruitment, it is now recognized as a much broader concept that includes attracting, retaining, developing and transitioning talented employees [6].

Talent management has been studied from different points of view (e.g., [7–9]). Among all the existing studies on characterizing talent management conceptually, Cappelli [2] has contributed a new perspective to talent management, based on the just-in-time manufacturing: a talent-on-demand framework. In this framework, forecasting product demand is comparable to forecasting talent needs. In particular, Cappelli [2] posits that estimating the cheapest and fastest ways to manufacture products is equivalent to cost-effectively developing talent for different sectors, including manufacturing, services, and technology. He also compares outsourcing to certain aspects of manufacturing processes, such as hiring external talent (on-time delivery) to effectively plan for succession events (delivery execution). The issues and challenges in managing the internal talent pipeline and advancing employees through jobs and experience developments are remarkably similar to the processes that occur in product movement through a supply chain. In both cases, companies benefit substantially from reducing bottlenecks that block advancement, speeding up processing times and improving forecasts to avoid mismatches. The benefits of Cappelli's new metaphor for talent management have been well characterized in the literature (e.g., [10–11]). However, since humans are distinct and more complex to model than inanimate machine parts in traditional supply chain models, effective application of Cappelli's perspective requires new models and methods to capture and represent these distinct features and complexities. Concepts such as hiring, training, retaining, learning and acquiring new skills fundamentally influence the analytical problems relevant to human capital management. Accordingly, these concepts are critical to the successful management and planning of a human resource supply chain in practice.

*Literature Review*

The definition of talent management given in the prior section spans all processes that include attraction and selection of required talents outside the organization through the retention and development of current personnel inside the organization. In this section, we review prior *analytical* studies implemented in each phase of the talent management process.

*Personnel selection*: Selecting personnel satisfying some assessment criteria has attracted considerable attention from scholars. Various multiple-criteria decision making (MCDM) approaches with precise or fuzzy evaluation criteria constitute the majority of analytical tools in this category. For instance, Gungor et al. [12] proposed a personnel selection system based on a fuzzy analytic hierarchy process (FAHP) to evaluate the best personnel using ratings of both qualitative and quantitative criteria. Lin [13] dealt with the personnel selection problem by integrating analytic network process (ANP) and fuzzy data envelopment analysis (DEA) approach. Gibney and Shan [14] described the use of AHP in the dean selection process and compared the results against the Provost's final decision. Liao and Chang [15] used ANP in the Taiwanese hospital public relations personnel selection process. In addition, Kelemenis and Askounis [16] applied fuzzy technique for order preference by similarity to an ideal solution (TOPSIS) for selection of qualified human resources. Dursun and Karsak [17] developed a fuzzy TOPSIS method in which both linguistic and numerical assessment scales were used in personnel selection problem. Kelemenis et al. [18] incorporated the TOPSIS method with fuzzy logic in order to select the most appropriate managers. Using Karnik–Mendel algorithm, Sang et al. [19] proposed an analytical solution to the fuzzy TOPSIS method to obtain fuzzy relative closeness in personnel selection problems. Kabak et al. [20] proposed a fuzzy hybrid MCDM approach including fuzzy ANP, fuzzy TOPSIS and fuzzy ELECTRE techniques to solve personnel selection problem.

*Career development/personnel advancement*: The purpose of personnel development is to enrich employees' skills and knowledge in order to make them fit for future positions, emergent roles and responsibilities. There are fewer studies in this problem category, and managerial subjective models are generally favored. For example, Hedge et al. [21] described the design of a career development and advancement system to compute an advancement score based on achievement of predefined milestones in the US Navy. Baruch [22] developed a career development model entitled Career Active System Triad (CAST) after reconciling the interests and concerns of both employers and employees. Hamori and Koyuncu [23] addressed the relationship between international assignments, a method of developing global leader competencies, and career advancement. By combining individual- and organization-level variables, Garavan et al. [24] investigated the factors predicting the career progression of hotel managers working in international hotel chains in Ireland, Europe and Asia.

*Personnel attrition/Churn/turnover*: Personnel attrition is another specific part of the talent management process that has attracted attention from management researchers. This problem is typically addressed using managerial and conceptual approaches. For instance, Morrell et al. [25] developed a theoretical and heuristic model describing the relationship between organizational change and employee's turnover. Using structural equation modeling and survival analysis, Hom and Kinicki [26] investigated how job dissatisfaction would progress into turnover. Maertz et al. [27] examined mediated effects of perceived supervisor support and perceived organizational support on turnover cognitions, and their interactive effects on turnover behavior. Steel and Lounsbur [28] and Holtom et al. [29] provided interested readers with reviews of past papers under this category.

*Employees' retention*: Arlotto, Chick, and Gans [30] analyzed the hiring and retention of heterogeneous workers who learned over time. Kyndt et al. [31] found that learning and development of employees are important retention-supporting strategies. Ramlall [32] provided a review of motivation theories and explained how employee motivation affects employee retention within organizations. Whitt [33] developed a mathematical model to help analyze the benefit in contact-center performance obtained from increasing employee retention. Hausknecht et al. [34] proposed and tested a model of 12 content-related factors thought to be partially responsible for employees' decisions to stay with a particular employer.

Although there might exist other papers in each category, we observe that the papers reviewed above only focus on a particular part of the talent management process. Consequently, the solution obtained in one phase that satisfies a specific set of decision criteria does not seem to function well in another phase of talent management under our definition. Hence, the literature suffers from lack of an integrated approach optimizing the entire end-to-end process of talent management. To eliminate this shortcoming, we seek in the present paper to develop a novel mathematical programming model, which enables us to model simultaneously all essential phases of the talent management process. These phases include attracting, interviewing and offering candidates for an organizational role, as well as the downstream phases of promoting and turnover of employees once employed. The model development is grounded on the modern perspective to talent management described by Capelli [2]. In this regard, we propose a multi-period mixed integer nonlinear programming model using a supply–demand perspective, in which, to take into account the natural uncertainty in the talent management process, some uncertain parameters are assumed stochastic following arbitrary probability density functions (PDFs). After linearizing the proposed nonlinear programming model to overcome any difficulties in deriving a solution, we employ chance-constrained programming (CCP) to deal with the stochastic parameters. CCP theory enables us to derive deterministic equivalent to one stochastic constraint. However, we need to exploit Monte Carlo simulation in using CCP theory to treat stochasticity involved in another constraint. This approach will be presented in a manner that uniquely addresses the problem of end-to-end talent management differently from prior approaches in the literature.

The paper unfolds as follows: Section 2 introduces the fundamental theories of TOPSIS and CCP. In Section 3, we briefly describe the underlined problem in this paper. Section 4 develops a stochastic multi-period mixed integer nonlinear programming model for talent pipeline management from a supply–demand standpoint. Section 5 presents a solution procedure for the proposed programming problem based on linearization of nonlinear terms involved, CCP and stochastic Monte Carlo simulation. Section 6 describes a case study in which we have succeeded in implementing our developed model. Finally, Section 7 is devoted to concluding remarks.

## 2. Fundamentals

### 2.1. TOPSIS

Any organization tries to identify the potential candidates for its job opportunities using limited number of recruitment channels. These channels may include career fairs, social media, personal interviews and communications, advertisements on different websites, and other approaches. However, based on its former experiences in employment and the quality of applications received through specific recruitment channels, an organization tends to prioritize all available channels. In this case, it may hire more people from channels with a higher priority when selecting potential candidates. In order to weight the recruitment channels in the present study based on pre-determined evaluation criteria, we exploit the classic TOPSIS method as a principal technique in multi attribute decision-making (MADM). It is grounded on the idea that the best alternative should logically have the shortest distance from the positive ideal solution and the farthest distance from the negative ideal solution. If each local criterion is monotonically increasing or decreasing, then it is easy to define an ideal solution. The positive (negative) ideal solution is achieved from all best (worst) attainable values of local criteria.

Suppose $m$ alternatives $A_1, A_2,..., A_m$ are going to be evaluated based on $n$ criteria $C_1, C_2,..., C_n$ in a MADM problem. Each alternative is evaluated with respect to these $n$ criteria. In this case, the decision matrix $D = [\pi_{ij}]_{m \times n}$ is used to present all ratings assigned to alternatives, where $\pi_{ij}$ is the rating of alternative $A_i$ in terms of the criterion $C_j$. Let $\Omega = (\omega_1, \omega_2,..., \omega_n)$ be the vector of local criteria weights satisfying $\sum_{j=1}^{n} \omega_j = 1$. Decision maker (DM) in this study gives these weights subjectively.

The TOPSIS method consists of the following steps [35]:

I. In the first step, we normalize the decision matrix as:

$$\overline{\pi}_{ij} = \frac{\pi_{ij}}{\sqrt{\sum_{k=1}^{m} \pi_{kj}^2}} \quad ; i = 1,...,m, \ j = 1,...,n. \quad (1)$$

Then, we multiply the columns of normalized decision matrix by the associated weights:

$$v_{ij} = \omega_j \times \overline{\pi}_{ij} \quad ; i = 1,...,m, \ j = 1,...,n. \quad (2)$$

II. The positive ideal and negative ideal solutions are, respectively, determined as follows:

$$A^+ = \{v_1^+, v_2^+,..., v_n^+\} = \{(\max_i v_{ij} \mid j \in K_b)(\min_i v_{ij} \mid j \in K_c)\} \quad (3)$$

$$A^- = \{v_1^-, v_2^-,..., v_n^-\} = \{(\min_i v_{ij} \mid j \in K_b)(\max_i v_{ij} \mid j \in K_c)\} \quad (4)$$

where $K_b$ and $K_c$ are the sets of benefit (positive) and cost (negative) criteria, respectively. In fact, the positive ideal and negative ideal solutions are composed of all best and worst values attainable of the criteria.

III. In the next step, we obtain the distances of the available alternatives from the positive ideal, $S_i^+$, and negative ideal, $S_i^-$, solutions. These Euclidean distances for each alternative are calculated as follows, respectively:

$$S_i^+ = \sqrt{\sum_{j=1}^{n}(v_{ij} - v_j^+)^2} \quad ; i = 1,...,m. \tag{5}$$

$$S_i^- = \sqrt{\sum_{j=1}^{n}(v_{ij} - v_j^-)^2} \quad ; i = 1,...,m. \tag{6}$$

IV. The relative closeness to the ideal alternatives is calculated as:

$$RC_i = \frac{S_i^-}{S_i^- + S_i^+} \quad ; i = 1,...,m, \quad 0 \leq RC_i \leq 1 \tag{7}$$

V. All alternatives are ranked according to the relative closeness to the ideal alternatives, i.e., the greater the $RC_i$ is, the more preferable the alternative $A_i$ is.

### 2.2. Chance-Constrained Programming; the Basic Theory

There are many management situations in which decisions are made in uncertain environments. In an uncertain environment, we are not able to precisely define the required parameters and gather their respective data to characterize a problem. In a prevalent uncertain environment, the required parameters of an optimization problem follow probabilistic distribution functions. In such cases, these parameters are so-called stochastic parameters. In order to deal with stochastic parameters in an optimization problem, the literature has widely evolved into application of stochastic programming approaches. One of the most commonly used methods is chance-constrained programming (CCP) initially proposed by Charnes and Cooper [36]. In CCP modeling of a stochastic decision system, we suppose the constraints will hold at least $\alpha$ percent of time. Here, $\alpha$ denotes the confidence level provided by the decision maker as a proper safety margin.

We consider the following mathematical programming problem with stochastic parameters:

$$\begin{aligned} & \max \ f(\mathbf{x}) \\ & s.t. \\ & g_j(\mathbf{x}, \xi) \leq 0 \quad ; j = 1,...,p \end{aligned} \tag{8}$$

where $\mathbf{x}$ is an *n*-dimensional decision vector, $\xi$ is a stochastic vector and $g_j(\mathbf{x},\xi); j=1,...,p$ are stochastic constraint functions. Since stochastic parameters are involved in characterizing all given constraints, we cannot reasonably define the *maximization* term in the objective function and the direction ≤ in the given constraints. In order to treat this type of problem correctly, Liu [37] suggested using the following CCP model:

$$\max \ f(\mathbf{x}) \tag{9}$$
$$s.t.$$
$$\Pr(g_j(\mathbf{x},\xi) \leq 0) \geq \alpha_j \quad ; j = 1,...,p \tag{10}$$

where Pr(.) implies the probability of the event in (.), and $\alpha_j$ is a predetermined confidence level to constraint *j*.

The CCP approach can easily incorporate more analysis and subjective assessments rather than other stochastic programming routines. In contrast to the two-stage stochastic programming approach in which violation of the constraints is allowed but penalized through a penalty term in the objective function, the CCP approach maintains a high level of reliability by expressing a minimum requirement on the probability of satisfying constraints. In other words, the resulting decision ensures the probability of complying with constraints, i.e., the confidence level of being feasible [38]. Since our framework involves principal decisions on human resources in the talent management problem, any violation of the organization's constraints may have critical consequences in terms of organizational budgets and operational-ability. Furthermore, the solution of a problem obtained by CCP provides comprehensive information on the economical outcome as a function of the desired confidence level of satisfying constraints, which is crucial for decision-making [38]. Hence, the CCP approach appears to be a better choice to treat uncertainty associated with the talent management problem.

## 3. Problem Description

As mentioned previously, we investigate the talent management problem in the present paper using the supply–demand philosophy proposed by Capelli [2]. In a supply–demand representation of the talent management problem, an organization acquires its essential human resources from available recruitment channels in each time period, which supply the required talents. Capelli [2] explained how an organization would benefit from modeling its talent management problem using this supply–demand perspective. Figure 1 depicts the talent management configuration schematically from this viewpoint.

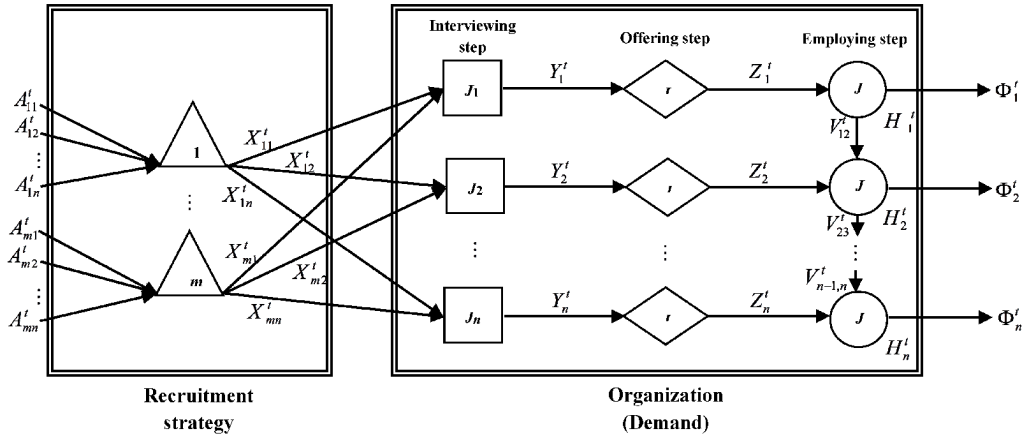

**Figure 1.** A schematic representation of the problem.

As shown in Figure 1, an organization continues accepting candidates for job position *j* in recruitment channel *i* in period *t* such that $A_{ij}^t$ candidates apply. Once the required number of applications is gathered for job position *j* in recruiting channel *i* during time period *t* (i.e., $A_{ij}^t$), the organization halts the application process and reviews the collected documents from the candidates. In this step, the organization invites $X_{ij}^t$ percent of all candidates with the highest qualifications for job position *j* applied via recruitment channel *i* to be interviewed in period *t*. Afterwards, it offers $Y_j^t$ percent of interviewed candidates (again with the highest qualifications for job position *j*) in time period *t*. Finally, $\tilde{b}_j$ percent of all offered candidates will accept their offers and start working in job position *j* in period *t*. To add additional complexity to this problem, in most organizations, there exist job seniorities based on organizational structure. In such structures, an individual can be promoted/demoted from a specific job position to another job position in each period based on his/her qualifications. The promotion/demotion is depicted in the far right tier in Figure 1 denoted by the advancement rate $V_{lj}^t$. Moreover, an employee working in position *j* may be dismissed by an organization for any reason or may quit of his/her own volition during period *t*. This fact is reflected by the attrition rate $\Phi_j^t$ in Figure 1. In this case, a vacant opportunity will become available for job position *j*. Ultimately, an organization will grow in size during the time horizon of the model. The growth of an organization is properly mirrored by a respective growth in its job positions, i.e., the more growth in job positions, the more growth the organization faces. We define the growth rate $G_j^t$ to denote the current condition. The key variables to be decided on in this study include the number of candidates for each job position applied through available recruitment channels, the number of required hires and the number of hired candidates along with operational rates (interview, offering, attrition, advancement and growth rates) in each time period in such a manner that the total average profit generated in an organization is maximized.

To understand and characterize the described problem better, we now present some remarks and relations below:

**Remark 1.** *If someone were promoted from job position l to job position j, then we would subtract 1 from the total number of employees in job position l and add 1 to the total number of employees in job position j.*

**Remark 2.** *Suppose someone is promoted from job position l to job position j. Based on the tacit relations in an organization, we can elicit the following relationships:*

(I)   # attritions for job position *j* = # current employees for job position *j* × attrition rate for job position *j*
(II)  # desired employees for job position *j* = # current employees for job position *j* × (1 + growth rate for job position *j*)
(III) # promotions *to* job position *j* = sum of (# current employees of any job position *l* × the advancement rate from job *l* to job *j*)
(IV)  # promotions *from* job position *j* = # current employees of job position *j* × sum of the advancement rates from job position *j* to any job position *k*
(V)   # gained/lost positions for job position *j* = # promotions *to* job position *j* − # attritions for job position *j* − # promotions *from* job position *j* = # promotions *to* job position *j* − # current employees for job position *j* × (attrition rate for job position *j* + sum of the advancement rates from job *j* to any job position *k*).

It is clear that if the value of the latter relation is positive, we have some employees gained for job position *j*; otherwise, some employees have been lost for that job position. Moreover, we have the following key relationships:

(VI)  # hires needed for job position *j* = # desired employees for job position *j* − # gained/lost positions for job position *j* − # current employees for job position *j* = # current employees of job position *j* × (growth rate for job position *j* + attrition rate for job position *j* + sum of the advancement rates *from* job *j* to any job position *k*) − # promotions *to* job position *j*
(VII) # of available employees in job position *j* at the end of a period = # current employees for job position *j* at the beginning of that period + # gained/lost positions for job position *j* + # of hired employees for job position *j*.

However, it should be noted here that an employee might not be promoted from job position *l* to job position *j* (or from job position *j* to job position *k*). For instance, an employee cannot be rationally promoted from the mechanical engineering position (*l*) to the market planning position (*j*) in an organization. To consider such a job relation (job seniority) based on organizational structure, a matrix called the transfer matrix is defined as $U = [u_{lj}]_{n \times n}$, where *n* is the total number of available job positions. Let us assume seniority of job position *j* is higher than that of job position *l*. In this case, if an employee can be promoted from the job position *l* to job position *j* (or degraded from job position *j* to job position *l*), then $u_{lj} = u_{jl} = 1$. On the other hand, if there is not such a job relation between job positions *j* and *l*, then $u_{lj} = u_{jl} = 0$.

## 4. Model Development

We propose a multi-period mixed integer nonlinear programming (MINLP) model with stochastic parameters to study the formerly described talent management problem. We introduce the required indices, parameters and decision variables as follows:

*Indices*

| | |
|---|---|
| $i$ | index of recruitment channels; $i = 1,...,m$. |
| $j$ | index of job positions; $j = 1,...,n$. |
| $t$ | index of time periods; $t = 1,...,T$. |

*Parameters*

| | |
|---|---|
| $\tilde{b}_j$ | the acceptance rate by candidates for job position *j*. |
| $\iota_j$ | the number of current employees for job position *j* at the beginning of planning horizon. |
| $e_j^t$ | the unit cost per hour of having extra hires for job position j in period t than it is required. |
| $\bar{e}_j^t$ | the unit cost per hour of having less hires for job position j in period t than it is required. |
| $\bar{\delta}_j^t$ | the maximum application capacity for job position j in period *t*. |
| $\bar{\beta}_j^t$ | the total interview rate for job position j in period *t*. |
| $\lambda_j^t$ | the maximum offering rate for job position *j* in period *t*. |
| $g_j^t$ | the maximum expected development for job position *j* in period *t*. |
| $\tilde{k}_j$ | the time (hour) required for considering initial documents and selecting candidates for being interviewed for job position *j*. |
| $\tilde{\tilde{k}}_j$ | the time (hour) required for analyzing the results of interviews for job position *j*. |
| $r_j^t$ | the unit revenue per hour yielded by each employee in job position j in period *t*. |
| $\psi_j^t$ | the unit salary per hour paid to each employee in job position j in period *t*. |
| $o_j^t$ | the unit cost of interview per hour for job position j in period *t*. |
| $\vartheta_j^t$ | the maximum number of employee changes for job position j in period *t*. |
| $w_i$ | the relative closeness of recruitment channel *i* derived from TOPSIS. |
| $\delta_i^t$ | the maximum application capacity of recruitment channel *i* in period *t*. |
| $\beta_i^t$ | the total amount of interview rate in recruitment channel *i* in period *t*. |
| $u_{jk}$ | the indicator parameter which is equal to 1 if someone can be promoted from job position *j* to job position *k* and vice versa, otherwise; is 0. |
| $RT^t$ | total recruitment process time (person-hour) in period *t*. |
| $\varepsilon$ | a very small positive number. |
| $M$ | a big positive number. |

It should be noted that the symbols with a tilde indicate parameters characterized by uncertainty. We characterize such parameters by appropriate probability distribution functions such as uniform, triangular, normal (Gaussian), etc. The main reason why we regard $\tilde{b}_j$, $\tilde{k}_j$ and $\tilde{\tilde{k}}_j$ as stochastic parameters is that their exact values are usually out of an organization's control. According to their definitions, $\tilde{k}_j$ and $\tilde{\tilde{k}}_j$ are directly affected by the number of applicants applied and their qualifications. Since an organization does not formerly know how many applicants are going to apply for a specific job position and what their qualifications would be, it cannot precisely determine the values of these time parameters. In addition, it is not clear in advance that how many offered applicants would accept their offerings. It is a subjective issue for applicants about which an organization does not rationally know. Hence, these points altogether prevent an organization to have any control on determining the exact values of $\tilde{b}_j$, $\tilde{k}_j$ and $\tilde{\tilde{k}}_j$ in data collection phase of the work.

*Decision variables*

$A_{ij}^t$  the number of candidates for job position *j* applied via recruitment channel *i* in period *t*.

$X_{ij}^t$  the interview rate per candidate for job position *j* in recruitment channel *i* in period *t*.

$Y_j^t$  the offering rate per candidate for job position *j* in period *t*.

$\Phi_j^t$  the attrition rate per employee for job position *j* in period *t*.

$V_{lj}^t$  the advancement rate per employee from job position *l* to job position *j* in period *t*.

$G_j^t$  the growth need by role for job position *j* in period *t*.

$S_j^t$  the number of employees for job position *j* at the end of period *t*.

$C_j^t$  the number of current employees for job position *j* at the beginning of period *t*.

$Z_j^t$  the number of candidates hired for job position *j* in period *t*.

$H_j^t$  the number of hires needed for job position *j* in period *t*.

$P_j^t$  the binary variable which is one (zero) if we have less (more) hires for job position *j* in period *t*.

Regarding these notes altogether and the problem statement, we construct a multi-period mixed integer nonlinear programming (MINLP) model as follows.

*4.1. Objective Function*

The objective function in the developed programming model seeks to maximize the total average profit per hour, which is defined as the difference between total average revenue per hour and total cost per hour. The total cost per hour is itself characterized as the summation of total recruitment process cost and total salary, both per hour. We formulize this objective function as follows:

$$\max f = \frac{1}{T}\left\{0.5\sum_{j=1}^{n}\sum_{t=1}^{T}r_j^t(S_j^t+C_j^t) - \{\sum_{j=1}^{n}\sum_{i=1}^{m}\sum_{t=1}^{T}o_j^t X_{ij}^t A_{ij}^t + \sum_{j=1}^{n}\sum_{t=1}^{T}[p_j^t \bar{e}_j^t(H_j^t-Z_j^t)+e_j^t(Z_j^t-H_j^t)(1-P_j^t)] + 0.5\sum_{j=1}^{n}\sum_{t=1}^{T}\psi_j^t(S_j^t+C_j^t)\}\right\} = $$
$$\frac{1}{T}\left\{0.5\sum_{j=1}^{n}\sum_{t=1}^{T}(r_j^t-\psi_j^t)(S_j^t+C_j^t) - \{\sum_{j=1}^{n}\sum_{i=1}^{m}\sum_{t=1}^{T}o_j^t X_{ij}^t A_{ij}^t + \sum_{j=1}^{n}\sum_{t=1}^{T}(P_j^t \bar{e}_j^t+(1-P_j^t)e_j^t)|H_j^t-Z_j^t|\}\right\}$$
(11)

Notably, the revenue generated by job position *j* at each period is calculated by multiplying its corresponding unit revenue per hour by the average number of available employees in that job position. Hence, $0.5\sum_{j=1}^{n}\sum_{t=1}^{T}r_j^t(S_j^t+C_j^t)$ denotes the total average revenue generated per hour over all job positions during all periods. We have assumed a linear relationship between the average number of available employees and the total average revenue they generate in an organization. Meanwhile, we could have calculated the total average revenue using a nonlinear regression model, in which the total average revenue yielded is estimated by the average number of all available employees. In addition, the total cost per hour represented in the accolade consists of interview costs, the cost of having both less and extra hires than required and total average salary paid to all job positions in all time periods.

*4.2. Constraints*

As stated before, one of the main aspects of talent management problem is to hire the right people at the right time. In order to address the issue of "time rightness", we set an upper bound to the time of recruitment process towards employing the "right people" in each period. When an organization needs someone to be employed, all necessary activities for employing him/her should be reasonably completed during a desirable time. In this regard, Constraint (12) bounds the total recruitment process time from applying to a specific job position to proposing offerings to the selected candidates in each period. This portion of recruitment process for the specific job position *j* contains (*a*) gathering, organizing and analyzing all documents related to applying candidates and inviting the top ones to interview; and (*b*) conducting interview, analyzing its outcome and recommending offers to finally selected candidates. The times associated with parts *a* and *b* equal $\tilde{k}_j \sum_{i=1}^m A_{ij}^t$ and $\tilde{\tilde{k}}_j \sum_{i=1}^m A_{ij}^t X_{ij}^t$ for job position *j* in period *t*, respectively. The summation of these two terms over all job positions yields to left hand side of Constraint (12).

$$\sum_{j=1}^n \sum_{i=1}^m (\tilde{k}_j + \tilde{\tilde{k}}_j X_{ij}^t) A_{ij}^t \leq RT^t \quad ; \forall t \tag{12}$$

Constraints (13) and (14) denote the relationships derived before in Parts VI and VII of Remark 2. Constraint (15) implies that the number of available employees at the end of some period *t* is naturally equal to the number of current employees at the beginning of next period, where $C_j^1 = \iota_j ; \forall j$.

$$H_j^t - C_j^t (G_j^t + \Phi_j^t + \sum_{\substack{k=1 \\ k \neq j}}^n u_{jk} V_{jk}^t) + \sum_{\substack{l=1 \\ l \neq j}}^n u_{lj} C_l^t V_{lj}^t = 0 \quad ; \forall j, t \tag{13}$$

$$S_j^t - Z_j^t - C_j^t (1 - \Phi_j^t - \sum_{\substack{k=1 \\ k \neq j}}^n u_{jk} V_{jk}^t) - \sum_{\substack{l=1 \\ l \neq j}}^n u_{lj} C_l^t V_{lj}^t = 0 \quad ; \forall j, t \tag{14}$$

$$S_j^{t-1} - C_j^t = 0 \quad ; \forall j, t \geq 2 \tag{15}$$

By consulting some HR managers in different industries, we noticed that it is not of interest for an organization to have its employees working in a particular job position change too much. In other words, an organization tends to retain its employees in a job position as long as possible. In this case, not only does an organization keep its uniformity in terms of the combination of its employees working in a job position, but also the employees enjoy job stability and are able to gain experience in a particular job position over time. To address this concern, we need to define an upper bound on the total number of changes in each job position, which is derived as the maximum of the total number of incoming employees to a job position and the total number of outgoing employees from that job position in each period. Hence, we have the following constraint:

$$\max(Z_j^t + \sum_{\substack{l=1 \\ l \neq j}}^n u_{lj} C_l^t V_{lj}^t , C_j^t (\Phi_j^t + \sum_{\substack{k=1 \\ k \neq j}}^n u_{jk} V_{jk}^t)) \leq \vartheta_j^t \quad ; \forall j, t \tag{16}$$

Furthermore, Constraint (17) is used to restrict the number of candidates hired for each job position in each period. In this regard, we know that $X_{ij}^t A_{ij}^t$ is the number of interviewed candidates for job position *j* in period *t*, who had applied via the recruitment channel *i*. Thus, $\sum_{i=1}^m X_{ij}^t A_{ij}^t$ gives the total number of interviewed candidates for job position *j*, from whom $Y_j^t$ percent are offered in period *t*. Now, $\tilde{b}_j$ percent of all offered candidates accept the offerings.

Hence, $\tilde{b}_j Y_j^t \sum_{i=1}^{m} X_{ij}^t A_{ij}^t$ yields the total number of candidates accepting the offerings for job position *j* in period *t*. However, for any reason, there might be situations in which an organization does not finally employ someone who has accepted an offering. This frequently happen especially when it does further considerations and withdraws its offering from a candidate. Considering such special cases, Constraint (17) implies the number of hired candidates in period *t* should be less than or equal to the number of accepted offerings.

$$Z_j^t - \tilde{b}_j Y_j^t \sum_{i=1}^{m} X_{ij}^t A_{ij}^t \leq 0 \quad ; \forall j, t \tag{17}$$

The condition in which we have either excessive hires or shortages is determined by Constraint (18). In this constraint, if the binary variable $P_j^t$ equals one (zero), then we have a hire shortage (excess) for the specific job position *j* in period *t*.

$$(Z_j^t - H_j^t) P_j^t + (H_j^t - Z_j^t)(1 - P_j^t) \leq 0 \quad ; \forall j, t \tag{18}$$

Certainly, each recruitment channel has a maximum capacity of application. For instance, we might be able to register a specific number of candidates for all job positions in a career fair in each period. To take into account this condition, Constraint (19) restricts the total number of candidates applying via recruitment channel *i* in period *t* to be less than or equal to $\delta_i^t$, where $w_i$ denotes the relative closeness of the recruitment channel *i* computed by TOPSIS. As it is observed, this constraint appropriately combines the priority and capacity of recruiting channel together.

$$\sum_{j=1}^{n} A_{ij}^t \leq w_i \delta_i^t \quad ; \forall i, t \tag{19}$$

Moreover, we impose a maximum capacity of application for a specific job position in all recruitment channels in each period. For instance, we might determine an upper bound for the number of candidates of job position financial clerk in all recruitment channels to be less than or equal to 15. On the other hand, when a job opportunity is announced, we would like to rationally have at least one candidate applied for it. Constraint (20) satisfies these conditions and sets upper and lower bounds on the total number of candidates for job position *j* applied via all recruitment channels in period *t*.

$$1 \leq \sum_{i=1}^{m} A_{ij}^t \leq \bar{\delta}_j^t \quad ; \forall j, t \tag{20}$$

In addition, Constraint (21) implies that the total amount of interview rates for all job positions in recruitment channel *i* in period *t* should be less than or equal to $\beta_i^t$. For example, an organization might desire having at most 25% of its total interviews from its website. Recall that this constraint appropriately uses the priority of recruitment channel *i* by combining the relative closeness $w_i$.

$$\sum_{j=1}^{n} X_{ij}^t \leq w_i \beta_i^t \quad ; \forall i, t \tag{21}$$

Constraint (22) limits the interview rate for specific job position *j* in all recruitment channels in period *t* to be less than or equal to $\bar{\beta}_j^t$. For instance, a company might want to have at most 30% of all candidates for financial clerk position interviewed during some period *t*.

$$\sum_{i=1}^{m} X_{ij}^t \leq \bar{\beta}_j^t \quad ; \forall j, t \tag{22}$$

Furthermore, we know that if no candidate for job position *j* in period *t* applies through recruitment channel *i*, i.e., $A_{ij}^t = 0$, then we will not have someone to be interviewed for that job position in recruitment channel *i*. In other words, the corresponding interview rate $X_{ij}^t$ will certainly equal zero in this case. Moreover, if job position *j* is not going to be offered to someone in period *t*, i.e., $Y_j^t = 0$, then it is economically reasonable that the interview rate for that job position will be zero, i.e., $\sum_{i=1}^m X_{ij}^t = 0$, which results in $X_{ij}^t = 0$. Hence, if either $A_{ij}^t = 0$ or $Y_j^t = 0$, then we must have $X_{ij}^t = 0$. On the other hand, if someone applies for job position *j* via recruitment channel *i* in period *t*, i.e., $A_{ij}^t > 0$, then we would like to interview with some applied candidates (at least one candidate). In this case, the respective interview rate should be greater than zero, i.e., $X_{ij}^t > 0$. The key Constraint (23) is to satisfy all aforementioned conditions, where *ε* and *M* are very small and big positive numbers, respectively.

$$\varepsilon A_{ij}^t \leq X_{ij}^t \leq M A_{ij}^t Y_j^t \quad ; \forall i,j,t \tag{23}$$

In addition, if the interview rate for job position *j* is equal to zero in period *t*, i.e., $\sum_{i=1}^m X_{ij}^t = 0$, then the respective offering rate for that job position should be zero, i.e., $Y_j^t = 0$, because there is no candidate available to be offered. On the other hand, if the interview rate for job position *j* is greater than zero in period *t*, i.e., $\sum_{i=1}^m X_{ij}^t > 0$, then we need to logically offer job position *j* to some applicants (at least one candidate), i.e., $Y_j^t > 0$. We imply these key conditions by Constraint (24).

$$\varepsilon \sum_{i=1}^m X_{ij}^t \leq Y_j^t \leq \min\{M \sum_{i=1}^m X_{ij}^t, \lambda_j^t\} \quad ; \forall j,t \tag{24}$$

Based on the definition of advancement rate, we impose Constraint (25) to restrict the summation of advancement rates from job position *j* to any job position *k* in each period to be less than or equal to one.

$$\sum_{k=1}^n u_{jk} V_{jk}^t \leq 1 \quad ; \forall t,j,k \neq j \tag{25}$$

Constraint (26) restricts the growth rate to lie within an interval whose upper bound reflects the maximum expected development for the respective job position in period *t*. Furthermore, it is known that all operational rates are between zero and one. The current requirement is satisfied by Constraint (27). In addition, according to the definitions of variables $A_{ij}^t, H_j^t, Z_j^t$ and $S_j^t$, $\forall i,j,t$ (Constraint (28)) satisfy their integrity conditions. Finally, following our discussion before, Constraint (29) requires that the decision variable $P_j^t$ be binary.

$$G_j^t \in [0, g_j^t] \quad ; \forall j,t \tag{26}$$

$$X_{ij}^t, \Phi_j^t, Y_j^t, V_{lj}^t \in [0,1] \quad ; \forall i,j,l,t \tag{27}$$

$$A_{ij}^t, H_j^t, Z_j^t, S_j^t, C_j^t : \text{integer} \quad ; \forall i,j,t \tag{28}$$

$$P_j^t \in \{0,1\} \quad ; \forall j,t \tag{29}$$

## 5. Solution Procedure

*5.1. Linearization*

The programming model developed in the preceding section is a multi-period mixed integer *nonlinear* problem, which might arise hardships to be solved. To conquer any of such difficulties, we linearize nonlinear terms that are present in the objective function and constraints by some mathematical programming techniques [39,40] and obtain an equivalent mixed integer *linear* programming (MILP) model.

The absolute term appeared in Objective Function (11), i.e., $|H_j^t - Z_j^t|$, can be eliminated from the model by substituting with a new integer variable $\gamma_j^t$ and adding two constraints as $\gamma_j^t \geq H_j^t - Z_j^t$ and $\gamma_j^t \geq Z_j^t - H_j^t$. Consequently, the integer decision variable $\gamma_j^t$ is multiplied by binary variable $p_j^t$ in the last term of the objective function. Now, the product $\gamma_j^t P_j^t$ denoted by $\overline{\gamma}_j^t$ can be linearized by adding the following four constraints: $\overline{\gamma}_j^t \leq M P_j^t$, $\overline{\gamma}_j^t \geq 0$, $\overline{\gamma}_j^t \leq \gamma_j^t$ and $\overline{\gamma}_j^t \geq \gamma_j^t - M(1 - P_j^t)$ for $\forall j,t$. Such a linearization is authentic because if $P_j^t = 0$, then the first two constraints imply that $\overline{\gamma}_j^t = 0$, and if $P_j^t = 1$, then the last two constraints gives $\overline{\gamma}_j^t = \gamma_j^t$.

Another nonlinear term in the objective function is $\sum_{j=1}^{n} \sum_{i=1}^{m} \sum_{t=1}^{T} o_j^t X_{ij}^t A_{ij}^t$ in which $X_{ij}^t$ and $A_{ij}^t$ are continuous and integer decision variables, respectively. Assuming $\hat{A}_{ij}^t = \min\{w_i \delta_i^t, \overline{\delta}_j^t\}$ with respect to Constraints (19) and (20), we determine a bound on $A_{ij}^t$ as $0 \leq A_{ij}^t \leq \hat{A}_{ij}^t$. In this case, the integer decision variable $A_{ij}^t$ is characterized as follows: $A_{ij}^t = \sum_{r_1=0}^{\overline{r}_1} 2^{r_1} (\overline{A}_{ij}^t)_{r_1}$ along with the constraint $A_{ij}^t \leq \hat{A}_{ij}^t$, where $\overline{r}_1 = \lfloor \log_2(\hat{A}_{ij}^t + 1) \rfloor$ and $(\overline{A}_{ij}^t)_{r_1} \in \{0,1\}$ for $\forall i,j,t,r_1$. Now, the abovementioned nonlinear term is rewritten as $\sum_{j=1}^{n} \sum_{i=1}^{m} \sum_{t=1}^{T} \sum_{r_1=0}^{\overline{r}_1} 2^{r_1} o_j^t X_{ij}^t (\overline{A}_{ij}^t)_{r_1}$, which contains multiplication of the continuous variable $X_{ij}^t$ by the binary variable $(\overline{A}_{ij}^t)_{r_1}$. To handle the current nonlinearity, we substitute $X_{ij}^t (\overline{A}_{ij}^t)_{r_1}$ by $(\overline{X}_{ij}^t)_{r_1}$, i.e., $(\overline{X}_{ij}^t)_{r_1} = X_{ij}^t (\overline{A}_{ij}^t)_{r_1}$, and add the following four sets of linear constraints: $(\overline{X}_{ij}^t)_{r_1} \leq (\overline{A}_{ij}^t)_{r_1}$, $(\overline{X}_{ij}^t)_{r_1} \geq 0$, $(\overline{X}_{ij}^t)_{r_1} \leq X_{ij}^t$ and $(\overline{X}_{ij}^t)_{r_1} \geq X_{ij}^t - (1 - (\overline{A}_{ij}^t)_{r_1})$ for $\forall i,j,t,r_1$. We also use the later procedure to eliminate the nonlinearity appeared in Constraint (12).

Using some arithmetic on Constraints (14)–(16), an upper bound for $C_j^t; \forall j,t \geq 2$ is derived as $\hat{C}_j^t = \max C_j^t = 2\vartheta_j^{t-1}; \forall j,t \geq 2$. Based on the current upper bound, the products $C_j^t G_j^t$, $C_j^t \Phi_j^t$ and $C_j^t V_{jk}^t$ in Constraints (13) and (14) can be linearized in the same way as described above. In this case, $C_j^t G_j^t$ is rewritten as $\sum_{r_2=0}^{\overline{r}_2} 2^{r_2} (\overline{C}_j^t)_{r_2} G_j^t$, where $(\overline{C}_j^t)_{r_2} \in \{0,1\}$ for $\forall j,r_2$ and $t \geq 2$, and $\overline{r}_2 = \lfloor \log_2(\hat{C}_j^t + 1) \rfloor$. Moreover, the following four sets of linear constraints are added: $(\overline{G}_j^t)_{r_2} \leq (\overline{C}_j^t)_{r_2} g_j^t$, $(\overline{G}_j^t)_{r_2} \geq 0$, $(\overline{G}_j^t)_{r_2} \leq G_j^t$ and $(\overline{G}_j^t)_{r_2} \geq G_j^t - (1 - (\overline{C}_j^t)_{r_2}) g_j^t$ where $(\overline{G}_j^t)_{r_2} = (\overline{C}_j^t)_{r_2} G_j^t$ for $\forall j,r_2$ and $t \geq 2$. To linearize the next quadratic term, we rewrite $C_j^t \Phi_j^t$ as $\sum_{r_2=0}^{\overline{r}_2} 2^{r_2} (\overline{C}_j^t)_{r_2} \Phi_j^t$ and then include the following four sets of linear constraints: $(\overline{\Phi}_j^t)_{r_2} \leq (\overline{C}_j^t)_{r_2}$, $(\overline{\Phi}_j^t)_{r_2} \geq 0$, $(\overline{\Phi}_j^t)_{r_2} \leq \Phi_j^t$ and $(\overline{\Phi}_j^t)_{r_2} \geq \Phi_j^t - (1 - (\overline{C}_j^t)_{r_2})$ where $(\overline{\Phi}_j^t)_{r_2} = (\overline{C}_j^t)_{r_2} \Phi_j^t$ for $\forall j,r_2$ and $t \geq 2$. Ultimately, we linearize the last quadratic term $c_j^t V_{jk}^t$ by rewriting it as $\sum_{r_2=0}^{\bar{r}_2} 2^{r_2} (\bar{C}_j^t)_{r_2} V_{jk}^t$ and imposing the following four sets of constraints: $(\bar{V}_{jk}^t)_{r_2} \leq (\bar{C}_j^t)_{r_2}$, $(\bar{V}_{jk}^t)_{r_2} \geq 0$, $(\bar{V}_{jk}^t)_{r_2} \leq V_{jk}^t$ and $(\bar{V}_{jk}^t)_{r_2} \geq V_{jk}^t - (1-(\bar{C}_j^t)_{r_2})$ where $(\bar{V}_{jk}^t)_{r_2} = V_{jk}^t (\bar{C}_j^t)_{r_2}$ for $\forall j, r_2$ and $t \geq 2$.

Constraint (16) can be transformed into two separate constraints as $Z_j^t + \sum_{l=1, l \neq j}^n u_{lj} C_l^t V_{lj}^t \leq \vartheta_j^t$ and $C_j^t (\Phi_j^t + \sum_{k=1, k \neq j}^n u_{jk} V_{jk}^t) \leq \vartheta_j^t$ for $\forall j, t$. We can linearize the product terms $C_l^t V_{lj}^t$, $C_j^t \Phi_j^t$ and $C_l^t V_{jk}^t$ in these separate constraints in same way as described before.

Constraint (17) seems more complicated to be linearized because it contains a product of one integer and two continuous decision variables. We first handle the product of two continuous decision variables, i.e., $X_{ij}^t Y_j^t$, and introduce two new continuous variables $\bar{y}_{ij}^t$ and $\bar{\bar{y}}_{ij}^t$, where $\bar{y}_{ij}^t = (X_{ij}^t + Y_j^t)/2$, $\bar{\bar{y}}_{ij}^t = (X_{ij}^t - Y_j^t)/2$, $0 \leq \bar{y}_{ij}^t \leq 1$ and $-1/2 \leq \bar{\bar{y}}_{ij}^t \leq 1/2$. In this case, we rewrite the product $X_{ij}^t Y_j^t$ as a separable function $(\bar{y}_{ij}^t)^2 - (\bar{\bar{y}}_{ij}^t)^2$. The quadratic terms of the current separable function are then piecewisely linearized using three breakpoints, which sounds a good approximation because the feasible intervals of the quadratic terms are not too wide. In this case, using three breakpoints of 0, 0.5 and 1, we obtain $(\bar{y}_{ij}^t)^2 = 0.25 \bar{\lambda}_{ij2}^t + \bar{\lambda}_{ij3}^t$ where $\bar{\lambda}_{ij1}^t + \bar{\lambda}_{ij2}^t + \bar{\lambda}_{ij3}^t = 1$, $\bar{\lambda}_{ij1}^t \leq \delta_{ij1}'^t$, $\bar{\lambda}_{ij2}^t \leq \delta_{ij1}'^t + \delta_{ij2}'^t$, $\bar{\lambda}_{ij3}^t \leq \delta_{ij2}'^t$, $\delta_{ij1}'^t + \delta_{ij2}'^t = 1$, $\bar{\lambda}_{ij1}^t, \bar{\lambda}_{ij2}^t, \bar{\lambda}_{ij3}^t \in [0,1]$ and $\delta_{ij1}'^t, \delta_{ij2}'^t \in [0,1]; \forall i, j, t$. In the same way, using three breakpoints of -1/2, 0 and 1/2, we have $(\bar{\bar{y}}_{ij}^t)^2 = 0.25 \bar{\bar{\lambda}}_{ij1}^t + 0.25 \bar{\bar{\lambda}}_{ij3}^t$ where $\bar{\bar{\lambda}}_{ij1}^t + \bar{\bar{\lambda}}_{ij2}^t + \bar{\bar{\lambda}}_{ij3}^t = 1$, $\bar{\bar{\lambda}}_{ij1}^t \leq \delta_{ij1}''^t$, $\bar{\bar{\lambda}}_{ij2}^t \leq \delta_{ij1}''^t + \delta_{ij2}''^t$, $\bar{\bar{\lambda}}_{ij3}^t \leq \delta_{ij2}''^t$, $\delta_{ij1}''^t + \delta_{ij2}''^t = 1$, $\bar{\bar{\lambda}}_{ij1}^t, \bar{\bar{\lambda}}_{ij2}^t, \bar{\bar{\lambda}}_{ij3}^t \in [0,1]$ and $\delta_{ij1}''^t, \delta_{ij2}''^t \in [0,1]; \forall i, j, t$. Consequently, the term $\sum_{i=1}^m Y_j^t X_{ij}^t A_{ij}^t$ in Constraint (17) is transformed into $\sum_{i=1}^m (0.25 \bar{\lambda}_{ij2}^t + \bar{\lambda}_{ij3}^t - 0.25 \bar{\bar{\lambda}}_{ij1}^t - 0.25 \bar{\bar{\lambda}}_{ij3}^t) A_{ij}^t$ with respect to all added constraints above. The later term contains products of continuous and integer variables and hence, as explained earlier, can be fully linearized. Assuming $\hat{\lambda}_{ij}^t = 0.25 \bar{\lambda}_{ij2}^t + \bar{\lambda}_{ij3}^t - 0.25 \bar{\bar{\lambda}}_{ij1}^t - 0.25 \bar{\bar{\lambda}}_{ij3}^t$ with lower and upper bounds of -0.25 and 1, respectively, we write $\hat{\lambda}_{ij}^t A_{ij}^t$ as $\sum_{r_1=0}^{\bar{r}_1} \hat{\lambda}_{ij}^t 2^{r_1} (\bar{A}_{ij}^t)_{r_1}$, then substitute $\hat{\lambda}_{ij}^t (\bar{A}_{ij}^t)_{r_1}$ with $(\breve{\lambda}_{ij}^t)_{r_1}$ and finally add the following four sets of linear constraints: $(\breve{\lambda}_{ij}^t)_{r_1} \leq (\bar{A}_{ij}^t)_{r_1}$, $(\breve{\lambda}_{ij}^t)_{r_1} \geq -0.25 (\bar{A}_{ij}^t)_{r_1}$, $(\breve{\lambda}_{ij}^t)_{r_1} \leq \hat{\lambda}_{ij}^t + 0.25((1-(\bar{A}_{ij}^t)_{r_1})$ and $(\breve{\lambda}_{ij}^t)_{r_1} \geq \hat{\lambda}_{ij}^t - (1-(\bar{A}_{ij}^t)_{r_1})$ for $\forall i, j, t, r_1$.

According to Constraint (16), an upper bound for $Z_j^t$ is obtained as $\vartheta_j^t$. Hence, in Constraint (18), we linearize the quadratic term $Z_j^t P_j^t$ by replacing it with $\bar{Z}_j^t$ where $\bar{Z}_j^t \leq \vartheta_j^t P_j^t$, $\bar{Z}_j^t \geq 0$, $\bar{Z}_j^t \leq Z_j^t$ and $\bar{Z}_j^t \geq Z_j^t - \vartheta_j^t (1-P_j^t)$. In addition, based on Constraints (13) and (16), and that $\max C_j^t = 2\vartheta_j^{t-1}; \forall j, t \geq 2$, an upper bound for $H_j^t$ can be obtained as $\hat{H}_j^t = \max H_j^t = 2\vartheta_j^{t-1} + \vartheta_j^t; \forall j, t$ where $\vartheta_j^0 = 0; \forall j$. Therefore, we substitute the product $H_j^t P_j^t$ in Constraint (18) with $\bar{H}_j^t$ where $\bar{H}_j^t \leq \hat{H}_j^t P_j^t$, $\bar{H}_j^t \geq 0$, $\bar{H}_j^t \leq H_j^t$ and $\bar{H}_j^t \geq H_j^t - \hat{H}_j^t (1-P_j^t)$.

We rewrite the quadratic term $A_{ij}^t Y_j^t$ in Constraint (23) as $\sum_{r_1=0}^{\bar{r}_1} 2^{r_1} Y_j^t (\overline{A}_{ij}^t)_{r_1}$, replace $Y_j^t (\overline{A}_{ij}^t)_{r_1}$ by $(\hat{y}_j^t)_{r_1}$ and add the following four sets of constraints: $(\hat{y}_j^t)_{r_1} \leq (\overline{A}_{ij}^t)_{r_1}$, $(\hat{y}_j^t)_{r_1} \leq Y_j^t$ and $(\hat{y}_j^t)_{r_1} \geq Y_j^t - (1 - (\overline{A}_{ij}^t)_{r_1})$ for $\forall i, j, t, r_1$.

Finally, we transform Constraint (24) to three usual linear constraints as: $\varepsilon \sum_{i=1}^{m} X_{ij}^t \leq Y_j^t$, $Y_j^t \leq M \sum_{i=1}^{m} X_{ij}^t$ and $Y_j^t \leq \lambda_j^t$ for $\forall j, t$.

### 5.2. Stochastic Constraints

According to the theory of CCP, it might be possible to convert the stochastic constraints to their deterministic equivalents for the predetermined confidence levels. Then, the deterministic programming problem is solved by typical solution approaches. Although this procedure seems relatively difficult to employ and only successful for special cases, we explain how it is effectively utilized along with Lemma 1 to deal with the stochastic Constraint (17) after being linearized. Nonetheless, we should use the stochastic simulation to encounter the stochastic Constraint (12), because the aforementioned conversion to obtain its deterministic equivalent after linearization is not easily applicable on it.

**Lemma 1.** *Assuming* $g(\mathbf{x}, \xi) = h(\mathbf{x}) - \xi$, *the deterministic equivalent of* $\Pr\{g(\mathbf{x}, \xi) \leq 0\} \geq \alpha$ *in CCP model is derived as* $h(\mathbf{x}) \leq \kappa_\alpha$, *where* $\kappa_\alpha = \varphi^{-1}(1-\alpha)$ *and* $\varphi^{-1}$ *is the inverse of cumulative distribution function* $\varphi(.)$ [37].

In the preceding section, we linearized the stochastic Constraint (17) as $Z_j^t - \tilde{b}_j \sum_{i=1}^{m} \sum_{r_1=0}^{\bar{r}_1} (\hat{\lambda}_{ij}^t)_{r_1} 2^{r_1} \leq 0$ ; $\forall j, t$ along with some additional constraints, which are avoided being rewritten here. Manipulating this linear constraint results in:

$$g(\mathbf{x}, \xi) = \frac{Z_j^t}{\sum_{i=1}^{m} \sum_{r_1=0}^{\bar{r}_1} (\hat{\lambda}_{ij}^t)_{r_1} 2^{r_1}} - \tilde{b}_j \leq 0 \quad ; \forall j, t \tag{30}$$

Regarding $Z_j^t / \sum_{i=1}^{m} \sum_{r_1=0}^{\bar{r}_1} (\hat{\lambda}_{ij}^t)_{r_1} 2^{r_1}$ and $\tilde{b}_j$ as $h(\mathbf{x})$ and $\xi$, respectively, we derive the following expression based on Constraint (10) in CCP theory and Lemma 1:

$$\Pr\{g(\mathbf{x}, \xi) \leq 0\} = \Pr\{\frac{Z_j^t}{\sum_{i=1}^{m} \sum_{r_1=0}^{\bar{r}_1} (\hat{\lambda}_{ij}^t)_{r_1} 2^{r_1}} \leq \tilde{b}_j\} \geq \alpha_1 \quad ; \forall j, t$$

$$= 1 - \varphi_{B_j}\left(\frac{Z_j^t}{\sum_{i=1}^{m} \sum_{r_1=0}^{\bar{r}_1} (\hat{\lambda}_{ij}^t)_{r_1} 2^{r_1}}\right) \geq \alpha_1 \quad ; \forall j, t \tag{31}$$

where $\varphi_{B_j}(.)$ is the cumulative distribution function of $\tilde{b}_j$. The latter expression can be further simplified as:

$$\frac{Z_j^t}{\sum_{i=1}^{m}\sum_{r_1=0}^{\bar{r}_1}(\hat{\lambda}_{ij}^t)_{r_1} 2^{r_1}} \leq \varphi_{B_j}^{-1}(1-\alpha_1) \quad ; \forall j, t \tag{32}$$

Now, we substitute the stochastic Constraint (17) with Constraint (32), which is its deterministic equivalent, in the proposed mathematical programming model. We could even more simplify Expression (32) providing that the stochastic variable $\tilde{b}_j$ followed a well-defined distribution function with a closed form cumulative distribution function. However, the exact distribution function of $\tilde{b}_j$ that can be exploited to derive its cumulative distribution function is available in rare situations. Therefore, we should infer in many cases the distribution function of $\tilde{b}_j$ from the available historical data. In such circumstances, the most fitted distribution function to the available data is estimated using statistical methods such as Chi-square or Kolmogorov–Smirnov test.

In addition, we linearized the stochastic Constraint (12) in the former section as $\sum_{j=1}^{n}\sum_{i=1}^{m}(A_{ij}^t\tilde{k}_j + \sum_{r_1=1}^{\bar{r}_1}\tilde{\tilde{k}}_j 2^{r_1}(\overline{X}_{ij}^t)_{r_1}) \leq RT^t; \forall t$ (again, we avoid citing the additional constraints here). Now, according to Constraint (10) in the CCP theory, we impose confidence level $\alpha_2$ on this constraint and rewrite it as follows:

$$\Pr\{\sum_{j=1}^{n}\sum_{i=1}^{m}(A_{ij}^t\tilde{k}_j + \sum_{r_1=1}^{\bar{r}_1}\tilde{\tilde{k}}_j 2^{r_1}(\overline{X}_{ij}^t)_{r_1}) \leq RT^t\} \geq \alpha_2 \tag{33}$$

For a very special case in which $\tilde{k}_j$ and $\tilde{\tilde{k}}_j$ are two *independent* stochastic parameters following *normal* distributions as $N(\mu_j, \sigma_j^2)$ and $N(\overline{\mu}_j, \overline{\sigma}_j^2)$ respectively, we are able to derive a deterministic equivalent to stochastic Constraint (33). In this case, we can prove that $\tau^t = \sum_{j=1}^{n}\sum_{i=1}^{m}(A_{ij}^t\tilde{k}_j + \sum_{r_1=1}^{\bar{r}_1}\tilde{\tilde{k}}_j 2^{r_1}(\overline{X}_{ij}^t)_{r_1})$ follows a normal distribution with mean and variance of $\mu_{\tau^t} = \sum_{j=1}^{n}\sum_{i=1}^{m}(A_{ij}^t\mu_j + \sum_{r_1=1}^{\bar{r}_1}\overline{\mu}_j 2^{r_1}(\overline{X}_{ij}^t)_{r_1})$ and $\sigma_{\tau^t}^2 = \sum_{j=1}^{n}\sum_{i=1}^{m}(A_{ij}^t\sigma_j^2 + \sum_{r_1=1}^{\bar{r}_1}\overline{\sigma}_j^2 2^{r_1}(\overline{X}_{ij}^t)_{r_1})$, respectively. Now, we can follow a somewhat similar argument in Lemma 1 to derive a deterministic equivalent to Constraint (33) as $RT^t \geq \mu_{\tau^t} + \sigma_{\tau^t}\kappa_{\alpha_2}$, where $\kappa_{\alpha_2} = \varphi_Z^{-1}(\alpha_2)$ and $\varphi_Z^{-1}(.)$ denotes the inverse of cumulative distribution function of a standard normal random variable. In this case, we substitute the stochastic Constraint (12) by the deterministic constraint $RT^t \geq \mu_{\tau^t} + \sigma_{\tau^t}\kappa_{\alpha_2}$ in the developed programming model.

However, the aforementioned very special case restricts us to a couple of strict assumptions: having independent random variables $\tilde{k}_j$ and $\tilde{\tilde{k}}_j$, and following the normal distributions. Based on the definitions of $\tilde{k}_j$ and $\tilde{\tilde{k}}_j$ in Section 4, we might not be able to infer in every situation that these two random variables are naturally independent to each other. More importantly, there might be lots of cases in which these two random variables do not follow a normal distribution. In such conditions, it might be extremely arduous to extract deterministic equivalent to the recently added stochastic Constraint (33). Instead, it is recommended applying Monte Carlo stochastic simulation. In the literature, there are several studies suggesting Monte Carlo stochastic simulation incorporated with genetic algorithm (GA) for solving chance constrained programming [41–43]. Meanwhile, we implement this kind of simulation in this paper using Lingo 16 optimization software [44], which provides us with a very suitable facility to handle the stochastic Constraint (33). In this regard, Appendix A, Part (a) presents the required built-in functions of Lingo 16 software to generally model a chance-constrained programming problem.

## 6. Empirical Study

In this section, we present how to use the proposed mathematical programming model and its solution procedure in a practical context. The model was validated with a logistics department of a Fortune 500 mining and construction corporation that designs, manufactures, markets and sells machinery and engines. The logistics department in this study was looking to hire candidates for five job positions including: coordinator, analyst, senior analyst, manager and senior manager ($j = 1,...,5$) in the course of three years ($T = 3$). In doing so, it usually uses three main recruitment channels to identify the potential candidates; $i = 1,2,3$. These channels consist of career fair, company website and social media. The career fairs are held regularly in a variety of universities across the United States at specific time slots. Furthermore, the company's website is always accessible for candidates to look for available job positions and apply to any of them. Moreover, this company also posts its job opportunities into some social media applications such as LinkedIn, CareerBuilder and so forth. The company expressed challenges in recruiting as it expanded its global presence and sought the advice of the research team, providing motivation for this study.

### 6.1. Data Acquisition

In our talent planning model, many parameters require that a company gather data from a number of sources. Data collection constitutes a fundamental step in implementing the proposed mathematical model in practice, so considerable attention to this step is needed to ensure data quality. The parameters described in Section 4 are categorized into two classifications based on their natures: deterministic and stochastic parameters. The most reasonable approach to quantifying deterministic parameters involved consulting with human resource experts to determine their most accurate estimates. On the other hand, relying on historical data and using statistical methods were considered appropriate for estimation of stochastic parameters, such as mean, variance, mode and so forth.

Table 1 presents the values of required parameters related to five job positions in the three time periods for this company. The number of current employees in this table, ($\iota_j$), reveals that a large logistics department in this company. It worth noting that $\iota_j$ is the number of employees presently working in the respective job position across the United States.

On the other hand, Table 2 shows the values of parameters associated with three recruitment channels over three time periods. As mentioned in Section 2.1, we defined criteria $C_1, C_2,...$ based on how recruitment channels are prioritized. The company identifies three criteria in this regard: the average of applicants' job experience ($C_1$), the average of applicants' requested salary per year ($C_2$) and the average score of applicants' university degree ($C_3$). From the company's standpoint, $C_1$ and $C_3$ are obviously positive criteria, while $C_2$ is a negative criterion. In order to obtain the average score of an applicant's university degree, we count the number of applicants with PhD, master, bachelor and associate degrees, who have already applied through each recruitment channel, and then assign a score of 9 to PhD degree, 7 to master degree, 5 to bachelor degree and 3 to associate degree. Then, we calculate $C_3$ by averaging the scores of all candidates. The quantities of each evaluating criterion for each recruitment channel are presented in Table 3. Following TOPSIS method, the company's experts assign equal weights of 1/3 to each criterion. In this case, the positive ideal and negative ideal solutions are determined as: $A^+ = \{0.2396, 0.1698, 0.2307\}$ and $A^- = \{0.1430, 0.2101, 0.1643\}$, respectively. Therefore, we have: $S_1^+ = 0.0966$, $S_2^+ = 0.0711$, $S_3^+ = 0.0889$, $S_1^- = 0.0777$, $S_2^- = 0.0977$, and $S_3^- = 0.0411$. Consequently, we derive the relative closeness (importance weight) of each recruitment channel as: 44.6%, 57.9% and 31.6%, respectively.

In the next step, we walk through estimating the stochastic parameters $\tilde{k}_j$, $\tilde{\bar{k}}_j$ and $\tilde{b}_j$; $j = 1,\ldots,5$. As mentioned before, the stochastic parameters are usually characterized using statistical methods such as Chi-square or the Kolmogorov–Smirnov test from the available historical data. We examined recorded data in the logistics department of this company from 1981 to 2013, a period during which the required data have been most reliably documented. Afterwards, some probability distribution functions are fitted to the available data sets of $\tilde{k}_j$, $\tilde{\bar{k}}_j$ and $\tilde{b}_j$; $j = 1,\ldots,5$ using Chi-square test. For instance, we fit a lognormal distribution function to $\tilde{k}_5$ for senior management job position and derived a probability plot as depicted in Figure 2. Since all data points are placed between 95% confidence interval and the *p*-value of the goodness of fit test is significantly greater than 0.05 as the critical value, we conclude based on Figure 2 that $\tilde{k}_5$ appropriately follows a lognormal distribution with location and scale parameters as 0.777 and 0.521, respectively.

**Table 1.** The required parameters associated with job positions in three time periods.

| Job Positions (*j*) | Time Period (*t*) | $\iota_j$ | $e_j^t$ | $\bar{e}_j^t$ | $\bar{\beta}_j^t$ | $\bar{\delta}_j^t$ | $\lambda_j^t$ | $r_j^t$ | $\psi_j^t$ |
|---|---|---|---|---|---|---|---|---|---|
| Coordinator | 1 | 125 | $5.90 | $47.73 | 0.2 | 1000 | 0.9 | $30.53 | $29.00 |
|  | 2 | - | $8.40 | $50.73 | 0.2 | 1000 | 0.9 | $32.20 | $30.68 |
|  | 3 | - | $9.40 | $53.78 | 0.2 | 1100 | 0.9 | $34.35 | $31.28 |
| Analyst | 1 | 96 | $7.40 | $49.50 | 0.2 | 700 | 0.9 | $37.80 | $36.25 |
|  | 2 | - | $10.40 | $52.50 | 0.2 | 700 | 0.9 | $39.48 | $37.93 |
|  | 3 | - | $11.15 | $54.22 | 0.2 | 770 | 0.9 | $41.63 | $38.53 |
| Senior analyst | 1 | 43 | $10.68 | $60.75 | 0.2 | 450 | 0.66 | $48.05 | $46.98 |
|  | 2 | - | $14.68 | $64.25 | 0.2 | 450 | 0.66 | $49.73 | $48.66 |
|  | 3 | - | $16.53 | $67.30 | 0.2 | 495 | 0.66 | $51.88 | $50.40 |
| Manager | 1 | 16 | $24.81 | $142.6 | 0.25 | 100 | 0.5 | $74.10 | $72.40 |
|  | 2 | - | $29.93 | $146.1 | 0.25 | 100 | 0.5 | $75.78 | $74.08 |
|  | 3 | - | $31.08 | $149.2 | 0.25 | 110 | 0.5 | $77.43 | $75.33 |
| Senior manager | 1 | 6 | $47.90 | $245.9 | 0.33 | 30 | 0.4 | $96.18 | $90.62 |
|  | 2 | - | $50.60 | $252.9 | 0.33 | 30 | 0.4 | $97.85 | $92.30 |
|  | 3 | - | $52.18 | $256.0 | 0.33 | 33 | 0.4 | $98.55 | $93.79 |

**Table 2.** The required parameters related to recruitment channels.

| Recruitment Channels (*i*) | Time Period (*t*) | $\beta_i^t$ | $\delta_i^t$ |
|---|---|---|---|
| Career fair | 1 | 0.65 | 1000 |
|  | 2 | 0.65 | 1080 |
|  | 3 | 0.65 | 1166 |
| Company website | 1 | 0.7 | 1000 |
|  | 2 | 0.6 | 1080 |
|  | 3 | 0.5 | 1166 |
| Social media | 1 | 0.8 | 500 |
|  | 2 | 0.75 | 540 |
|  | 3 | 0.75 | 583 |

**Table 3.** Decision matrix for evaluating recruitment channels.

| Recruitment Channels | $C_1$ (Year) | $C_2$ ($) | $C_3$ |
|---|---|---|---|
| Career fair | 1.85 | $56,000 | 7.61 |
| Company website | 3.10 | $64,400 | 5.42 |
| Social media | 2.36 | $69,300 | 5.80 |

We pursue the same procedure to identify the probability distribution functions of all remained stochastic parameters, as given in Table 4:

**Table 4.** Characterization of PDFs of stochastic parameters.

| Job Positions | $\tilde{k}_j$ | $\tilde{\tilde{k}}_j$ | $\tilde{b}_j$ |
| --- | --- | --- | --- |
| Coordinator | exp (2.6570) | exp (1.2091) | uniform (0.06,1.00) |
| Analyst | exp (1.3422) | exp (0.8482) | uniform (0.16,0.87) |
| Senior analyst | exp (1.1328) | exp (0.7617) | uniform (0.42,0.82) |
| Manager | exp (0.9961) | exp (0.6957) | uniform (0.72,1.00) |
| Senior manager | Lognormal (0.777,0.521) | Lognormal (1.019,0.467) | uniform (0.83,1.00) |

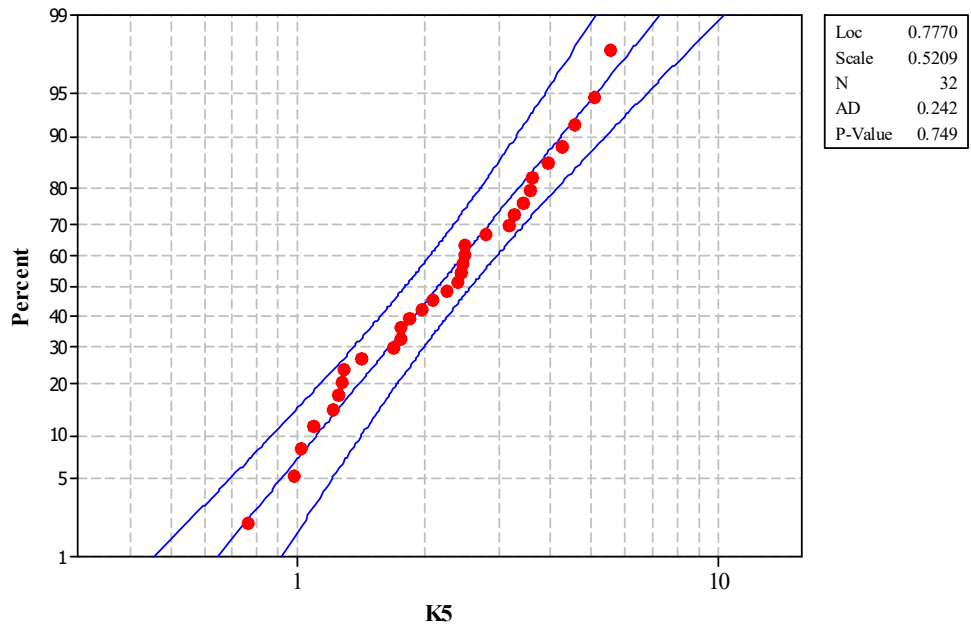

**Figure 2.** Lognormal probability plot of $\tilde{k}_5$.

As stated in Section 3, we sought to define the indicator parameter $u_{jk}$ to consider organizational structure in promoting an employee from job position *j* to *k* and vice versa. Table 5 presents the values of this indicator parameter for five available job positions. According to this table, an employee can be promoted to or demoted only from adjacent job positions, i.e., there is no skipping of levels in promotion or demotion processes in this department. For instance, since $u_{Analyst,Manager} = 0$, an employee cannot be promoted from analyst position to manager position in this department and vice versa, i.e., a manager is not demoted directly to an analyst position.

The next parameter identified in this study is the total number of employee changes for a job position in each period, i.e., $\vartheta_j^t$. In each period, the company sought to have employees changed in a job position by at most 50% of the respective number of employees in that position at the beginning of the period. In other words, $\vartheta_j^t = 0.5 C_j^t; \forall j,t$. At the end of data acquisition phase, we determined three remained parameters including $\varepsilon$, $M$ and $RT^t$ as 0.001, 10,000 and 480 person-hour, $t = 1,2,3$, respectively.

**Table 5.** Transfer matrix for five available job positions.

|  | Coordinator | Analyst | Senior Analyst | Manager | Senior Manager |
|---|---|---|---|---|---|
| Coordinator | - | 1 | 0 | 0 | 0 |
| Analyst | 1 | - | 1 | 0 | 0 |
| Senior analyst | 0 | 1 | - | 1 | 0 |
| Manager | 0 | 0 | 1 | - | 1 |
| Senior manager | 0 | 0 | 0 | 1 | - |

*6.2. Results Analysis*

Up to this step, we succeeded in gathering/defining all required parameters for the developed mathematical programming model. The next step was to solve the model by taking advantage of the stochastic programming facility in Lingo 16 optimization software [44] on a dual-core 2.5 GHz computer with 6 GB RAM. Table 6 presents the dimension of both original MINLP and linearized MILP models with the solution time of the latter model. According to this table, it appears that the dimension of MILP model is greater than that of original MINLP model. However, since the derived MILP model can effectively be solved using the available algorithms, its solution time does not sound irritating.

**Table 6.** Dimension of original MINLP and linearized MILP models.

|  | Original MINLP | Linearized MILP |
|---|---|---|
| # of decision variables | 103 | 839 |
| # of constraints | 113 | 1153 |
| Solution time (sec) | - | 89.3 |

Appendix A, Part (b) represents all the required commands of Lingo 16 that are associated with modeling stochasticity of the proposed mathematical programming problem. In this case, Tables 7 and 8 present the optimal values for $A_{ij}^t$ and $X_{ij}^t$, respectively:

**Table 7.** The optimal solutions of $A_{ij}^t$ in three time periods.

| Job Positions | Recruitment Channels | | | | | | | | |
|---|---|---|---|---|---|---|---|---|---|
|  | Career Fair | | | Company Website | | | Social Media | | |
|  | 1 | 2 | 3 | 1 | 2 | 3 | 1 | 2 | 3 |
| Coordinator | 0 | 92 | 78 | 82 | 0 | 0 | 0 | 0 | 0 |
| Analyst | 0 | 88 | 66 | 92 | 0 | 0 | 0 | 0 | 0 |
| Senior analyst | 15 | 0 | 0 | 0 | 27 | 21 | 0 | 0 | 0 |
| Manager | 0 | 0 | 0 | 0 | 49 | 56 | 27 | 0 | 0 |
| Senior manager | 0 | 0 | 0 | 0 | 0 | 0 | 18 | 14 | 33 |

**Table 8.** The optimal solutions of $X_{ij}^t$ in three time periods.

| Job Positions | Recruitment Channels | | | | | | | | |
|---|---|---|---|---|---|---|---|---|---|
|  | Career Fair | | | Company Website | | | Social Media | | |
|  | 1 | 2 | 3 | 1 | 2 | 3 | 1 | 2 | 3 |
| Coordinator | 0 | 0.092 | 0.079 | 0.198 | 0 | 0 | 0 | 0 | 0 |
| Analyst | 0 | 0.198 | 0.067 | 0.194 | 0 | 0 | 0 | 0 | 0 |
| Senior analyst | 0.187 | 0 | 0 | 0 | 0.104 | 0.053 | 0 | 0 | 0 |
| Manager | 0 | 0 | 0 | 0 | 0.204 | 0.223 | 0.092 | 0 | 0 |
| Senior manager | 0 | 0 | 0 | 0 | 0 | 0 | 0.158 | 0.217 | 0.237 |

For instance, in period one, the company collects 82 applications of candidates for coordinator position via its website. After that, it will not accept any application for this job position in its website and will remove the associated job opportunity announcement from the website. Hence, in period one, the total number of candidates for job positions coordinator, analyst, senior analyst, manager and senior manager invited to be interviewed is computed as 16.24 (=0.198 × 82), 17.85 (=0.194 × 92), 2.81 (=0.187 × 15), 2.48 (=0.092 × 27) and 2.84 (=0.158 × 18), respectively. Obviously, applicants with the highest qualifications are selected in this step to be interviewed. We can round these quantities as 16, 18, 3, 2 and 3 to be meaningful in practice. In addition, Table 9 gives the optimal solutions of the operational rates in conjunction with $Z_j^t$, $H_j^t$ and $S_j^t$:

**Table 9.** The optimal solutions of operational rates, $Z_j^t$, $H_j^t$ and $S_j^t$ in three time periods

| Job Positions | $\Phi_j^t$ | | | $G_j^t$ | | | $Y_j^t$ | | | $Z_j^t$ | | | $H_j^t$ | | | $S_j^t$ | | |
|---|---|---|---|---|---|---|---|---|---|---|---|---|---|---|---|---|---|---|
| | 1 | 2 | 3 | 1 | 2 | 3 | 1 | 2 | 3 | 1 | 2 | 3 | 1 | 2 | 3 | 1 | 2 | 3 |
| Coordinator | 0 | 0 | 0 | 0 | 0 | 0 | 0.90 | 0.69 | 0.90 | 5 | 2 | 1 | 5 | 2 | 1 | 125 | 125 | 125 |
| Analyst | 0 | 0 | 0 | 0.073 | 0.049 | 0.046 | 0.90 | 0.90 | 0.90 | 6 | 5 | 1 | 6 | 5 | 1 | 103 | 108 | 113 |
| Senior analyst | 0 | 0 | 0 | 0 | 0 | 0 | 0.66 | 0.66 | 0.66 | 1 | 1 | 1 | 1 | 1 | 1 | 43 | 43 | 43 |
| Manager | 0 | 0 | 0 | 0.313 | 0.286 | 0 | 0.50 | 0.50 | 0.50 | 1 | 4 | 5 | 1 | 4 | 5 | 21 | 27 | 27 |
| Senior manager | 0 | 0 | 0 | 0.333 | 0.250 | 0.500 | 0.40 | 0.37 | 0.40 | 1 | 1 | 2 | 1 | 1 | 2 | 8 | 10 | 15 |

The retention rate of one (attrition rate of zero) in period one given in this table implies that the company should retain all its current employees. Furthermore, according to the optimal values of offering rates in period one, the company is supposed to offer the job positions to 14.4 (=0.9 × 16), 16.2 (=0.9 × 18), 1.98 (=0.66 × 3), 1 (=0.5 × 2) and 1.2 (=0.4 × 3) interviewed candidates with the highest qualifications, respectively. Again, we may round these values as 14, 16, 2, 1 and 1 in practice. The offered candidates will respond to their offerings based on their respective stochastic acceptance rates given as $\tilde{b}_j$ in Table 4. The number of candidates hired in period one for each abovementioned job position is 5, 6, 1, 1 and 1, presented as $Z_j^1$. This implies the actual acceptance rates of offered candidates for each job position in period one are equal to 0.3571 (=5/14), 0.375 (=6/16), 0.5 (=1/2), 0.1 (=1/1) and 0.1 (=1/1), respectively. We observe that the actual acceptance rate for each job position in period one drops within its corresponding range of uniform distribution function given in Table 4.

Moreover, the optimal values of advancement rates are presented in Table 10:

**Table 10.** The optimal solutions of advancement rates in three time periods.

| Job Positions | Coordinator | | | Analyst | | | Senior Analyst | | | Manager | | | Senior Manager | | |
|---|---|---|---|---|---|---|---|---|---|---|---|---|---|---|---|
| | 1 | 2 | 3 | 1 | 2 | 3 | 1 | 2 | 3 | 1 | 2 | 3 | 1 | 2 | 3 |
| Coordinator | 0 | 0 | 0 | 0.040 | 0.016 | 0.008 | 0 | 0 | 0 | 0 | 0 | 0 | 0 | 0 | 0 |
| Analyst | 0 | 0 | 0 | 0 | 0 | 0 | 0.042 | 0.019 | 0 | 0 | 0 | 0 | 0 | 0 | 0 |
| Senior analyst | 0 | 0 | 0 | 0 | 0 | 0.069 | 0 | 0 | 0 | 0.116 | 0.07 | 0 | 0 | 0 | 0 |
| Manager | 0 | 0 | 0 | 0 | 0 | 0 | 0 | 0 | 0.074 | 0 | 0 | 0 | 0.063 | 0.047 | 0.111 |
| Senior manager | 0 | 0 | 0 | 0 | 0 | 0 | 0 | 0 | 0 | 0 | 0 | 0 | 0 | 0 | 0 |

Using the rates given in the current table, we are able to calculate the number of employees to be promoted to/demoted from each job position in each time period. For instance, in period one, since $V_{12}^1 = 0.040$, $V_{23}^1 = 0.042$, $V_{34}^1 = 0.116$, and $V_{45}^1 = 0.063$, the company should promote five (= 0.04 × 125) coordinators to analysts, 4.03 (= 0.042 × 96) analysts to senior analysts, 4.99 (= 0.116 × 43) senior analysts to managers and 1.01 (= 0.063 × 16) managers to senior managers. Reasonably, after rounding to the nearest integers, we select the employees with the highest qualifications for promotion. As observed, no degradation happens in this department during the first time period.

We interpreted the obtained results only for the first time period. Following the same procedure, one is able to understand the given results in the remained periods.

Table 9 also represents the number of required hires for each job position in each period, i.e., $H_j^t$. Comparing these values with the number of hired candidates in the respective time period, i.e., $Z_j^t$, we understand in the optimal condition that the company should employ the exact number of candidates it needs, i.e., no excess or shortage in employment is observed in this study. This observation refers to the optimal equilibrium between recruiting on the external labor market and the training and development of internal candidates.

Under the current optimal condition, the employees in the logistics department will yield $461.86, $490.29 and $904.61 profit per hour in consecutive time periods with an average of $618.92. The reason why we observe a significant increase in period three is the increase in net profit per hour, which is defined as the difference between generated revenue and paid salary per hour, yielded by a coordinator and an analyst in this period. According to Table 1, the net profits generated by a coordinator and an analyst in the third period are equal to $3.07 (=$34.35 − $31.28) and $3.10 (=$41.63 − $38.53) per hour, respectively, which both are almost twice as much the respective profits in the preceding periods. Furthermore, assuming 250 business days per year each with eight working hours, the logistics department in this company will annually generate $923,720, $980,580 and $1,809,220 in consecutive periods and consequently $3,713,520 in the whole time horizon.

Exploring the results obtained displays how the proposed approach in this paper is capable of managing the required talents in an organization. The right number of candidates for any job position from any recruiting channel in any time period was determined by introducing the decision variable $A_{ij}^t$; $\forall i, j, t$. Then, highly qualified persons were determined as the right candidates (employees) to be interviewed and offered (promoted). If necessary, inside an organization, the right employees to fire were selected among the less qualified persons in a particular job position. Moreover, the right place, i.e., job position, to employ the required candidates in each time period was also determined by introducing the respective decision variable $H_j^t$; $\forall j, t$. The issue of time rightness in hiring the required employees was addressed by setting an upper bound for the total time of recruiting activities, which guarantees they would be certainly hired during a desired time in each period. The results of the model were presented to the senior logistics executive at the company, as well as the human resources team responsible for recruiting and placement. The team agreed the model provided them with a much more concrete set of estimates for planning their recruiting strategies, and subsequently used to plan university events and website placement advertisements in the following year. Managers also commented that the model was relatively easy to use, as an Excel spreadsheet allowed easy input of parameters and reporting of results.

*6.3. Stochastic Analysis*

As stated in Appendix A, Part (b), we set the stochastic programming module of Lingo 16 optimization software to generate 60 random numbers of each stochastic parameter. Figure 3 illustrates the histograms of generated random numbers. In addition, the sample mean and standard deviation of generated random numbers associated with each stochastic parameter along with the *p*-values of goodness of fit tests are presented in Table 11. Examining Figure 3 and the statistics given in Table 11, we comprehend that the generated random numbers follow appropriately their respective theoretical distribution functions, as presented in Table 4.

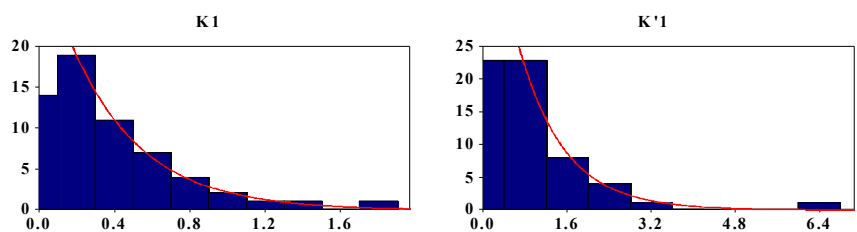

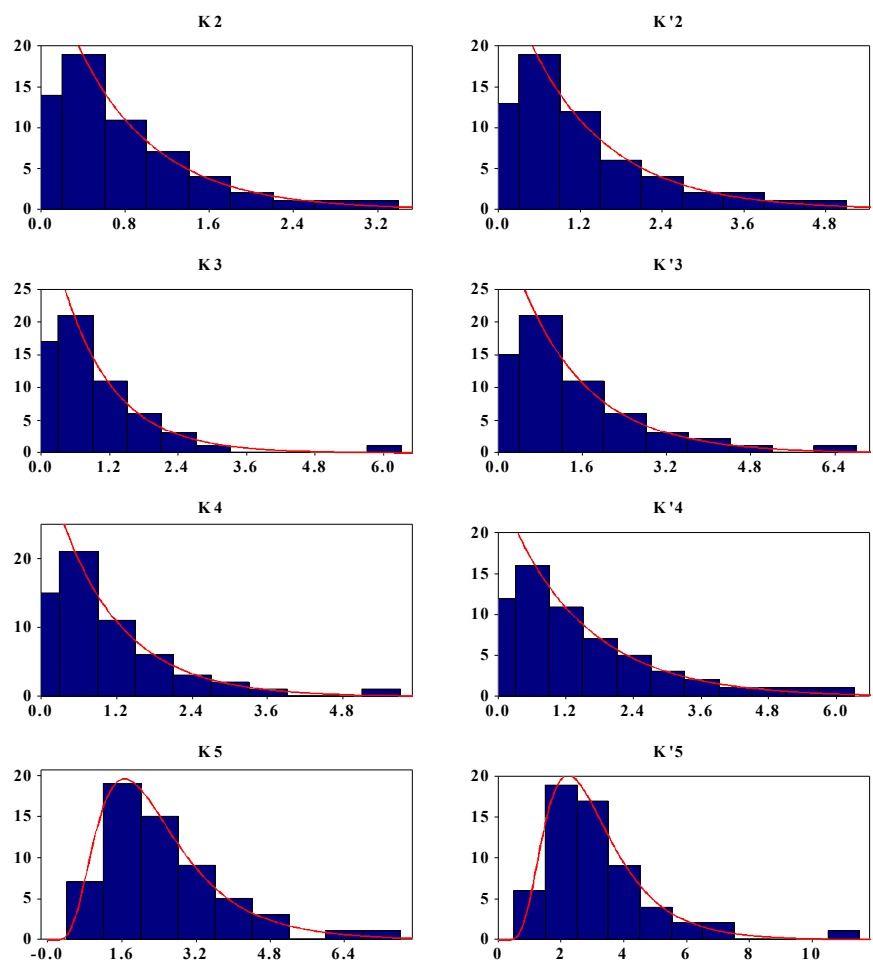

**Figure 3.** The histograms of generated random numbers by Lingo.

**Table 11.** The characteristics of generated random numbers by Lingo.

| Job Positions | $\tilde{k}_j$ | | | $\tilde{\tilde{k}}_j$ | | |
|---|---|---|---|---|---|---|
| | Sample Mean | Sample StdDev | *p*-Value | Sample Mean | Sample StdDev | *p*-Value |
| Coordinator | 0.3725 | 0.3597 | 1.00 | 0.8718 | 1.0136 | 0.997 |
| Analyst | 0.7382 | 0.7177 | 1.00 | 1.1606 | 1.0972 | 1.00 |
| Senior analyst | 0.9083 | 1.0021 | 1.00 | 1.3012 | 1.2729 | 1.00 |
| Manager | 1.0033 | 1.0017 | 1.00 | 1.4028 | 1.3258 | 1.00 |
| Senior manager | 2.4700 | 1.2948 | 1.00 | 3.1183 | 1.6386 | 1.00 |

In the theory of stochastic programming, there is a measure used for evaluating the importance of randomness involved in the problem: Expected Value of Perfect Information (EVPI). EVPI indicates the maximum amount a decision maker would be ready to pay in return for complete information about the future [45]. Although the concept of EVPI was initially developed in the context of decision analysis, it can be derived in the stochastic programming setting as follows [45]: for each set of realizations of stochastic parameters in the problem, we solve the stochastic program separately to obtain the corresponding objective function values. In this case, if we realize the stochastic parameters $n_s$ times (as the sample size), then $n_s$ different values of objective function will be attained, say $\zeta_l$; $l = 1,...,n_s$, which is known as the *distribution problem*. The expected value of the current quantities, i.e., $\bar{\zeta} = \sum_{l=1}^{n_s} \zeta_l / n_s$, is called the *wait-and-see solution* (WS). If value of the objective function in the original stochastic program is denoted by *VRP*, then $EVPI = VRP - WS$.

Figure 4 depicts values of the distribution problem in our case study, i.e., $\zeta_l; l = 1,\ldots,60$. The red and black horizontal lines in this figure demonstrate expected value of the distribution problem, which is equal to $615.33 per hour, and *VRP*, respectively. Hence, we conclude that the company will be willing to pay $3.59 (=$618.92 − $615.33) per hour to achieve the perfect information for future.

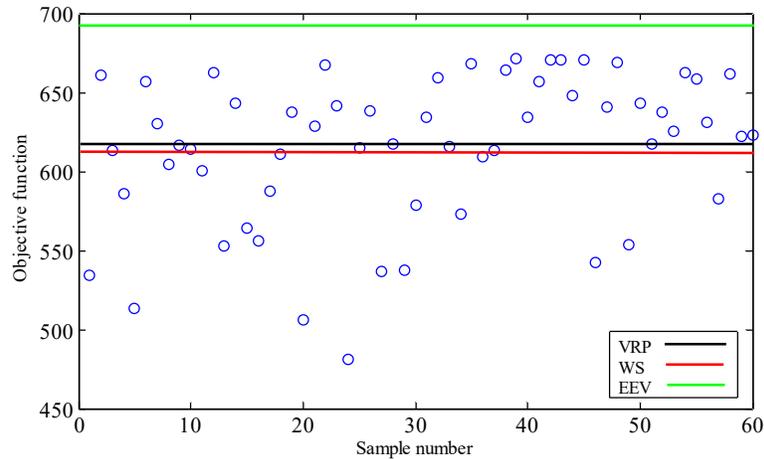

**Figure 4.** The values of distribution problem with WS and VRP.

Wait-and-see solution (*WS*) provides a lower bound for objective value of the original stochastic problem (*VRP*). We are also able to achieve an upper bound for it by substituting the stochastic parameters with their corresponding mean values. The obtained problem, called *mean value problem* (*EEV*), is of special interest in practice particularly once perfect information is just not available at any price [45]. If we solve the proposed model using the mean values of stochastic parameters given in Table 4, we obtain an upper bound of $693.40 for this case study. Hence, we derive the upper and lower bounds on value of objective function in the stochastic program as $693.40 and $615.33 per hour, respectively.

*6.4. Sensitivity Analysis*

As observed before, a long list of parameters needs to be gathered/defined in developing the proposed mathematical model. Now, it is of our interest here to investigate more profoundly the effect of changing in some substantial parameters on the given results, particularly the obtained value of objective function.

Two parameters that are set up intuitively in this case study are the confidence levels $\alpha_1$ and $\alpha_2$, which were initially determined to be equal to 0.7 and 0.95, respectively. Figure 5 depicts how the value of objective function changes if $\alpha_1$ and $\alpha_2$ individually vary within their corresponding feasible ranges. In Figure 5a, we observe if the confidence level $\alpha_1$ increases from 0.3 to 0.9, then the total average profit generated in this department decreases. The reasoning behind this observation is that if $\alpha_1$ increases, then the left hand side of Constraint (32) that is the inverse of cumulative distribution function of some uniform distributions becomes smaller. As a result, the value of maximization objective function does not become better, as it is decreasing here. Furthermore, under any specific value for $\alpha_1$, the generated profit in consecutive time periods is increasing, which is mainly because of increasing net profits generated by employees over successive time periods.

In addition, once the confidence level $\alpha_2$ increases from 0.1 to 0.95, the value of objective function decreases as well. It follows because, according to Constraint (33), increasing $\alpha_2$ makes the feasible region of problem more constrained, which impedes the objective function getting better. Again, due to increasing net profit yielded by employees over consecutive periods, we observe that, for each particular value of $\alpha_2$, the profit generated in this department also increases over successive time periods, as presented in Figure 5b.

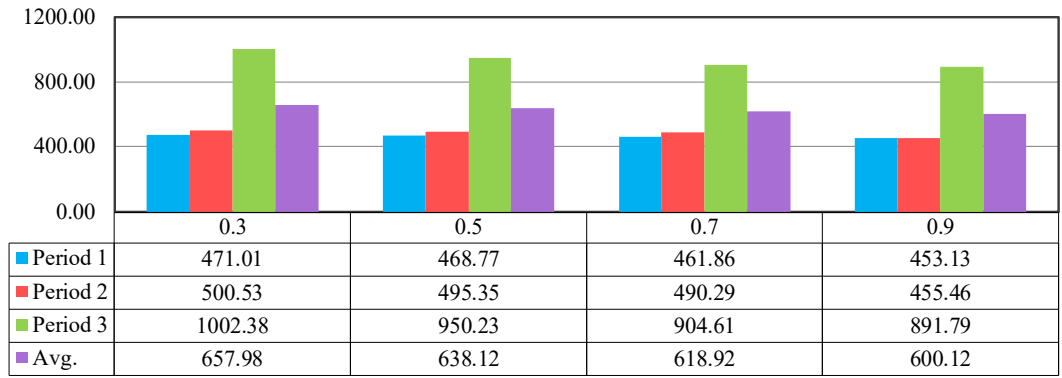

(**a**)

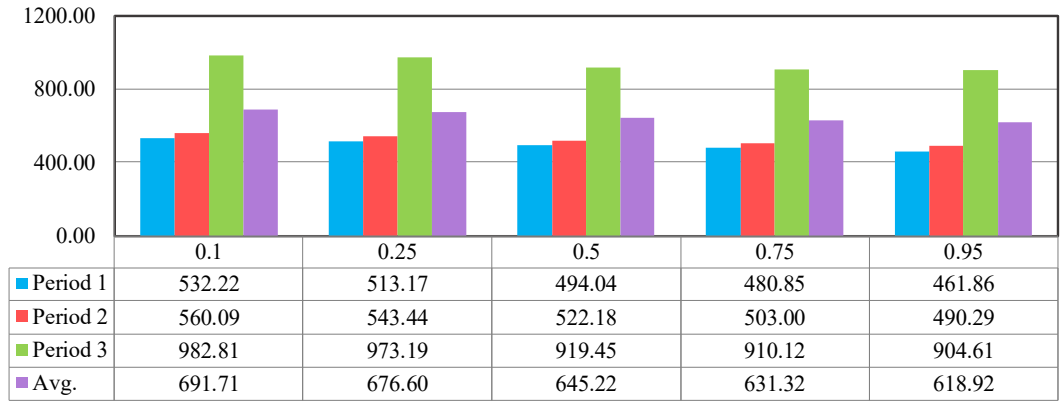

(**b**)

**Figure 5.** Sensitivity analysis on: $\alpha_1$ (**a**); and $\alpha_2$ (**b**).

As previously stated, we rely on historical data in conjunction with Chi-square test to derive the PDF of stochastic parameters including $\tilde{k}_j$, $\tilde{\tilde{k}}_j$ and $\tilde{b}_j$; $j=1,...,5$. Meanwhile, it is likely the future behavior of these uncertain parameters be somewhat different from what we characterize in Table 4. In order to encounter such volatility, having the confidence levels $\alpha_1$ and $\alpha_2$ fixed at their initial values, now we explore the influence of changes in the mean parameters of the derived PDF's presented in Table 4 on the generated total average profit in this department. To do so, we assume the mean of the stochastic parameters above can decrease and increase at most by 50%. In this case, Figure 6 depicts how sensitive the value of objective function is once the changes in the mean of time parameters $\tilde{k}_j$ and $\tilde{\tilde{k}}_j$, and acceptance rates $\tilde{b}_j$; $j=1,...,5$ occur.

We observe that if the mean of either time parameter decreases (increases), then the total average profit yielded in this department increases (decreases). We reason that once the mean of a time parameter (either $\tilde{k}_j$ or $\tilde{\tilde{k}}_j$) reduces (increases), more (less) applicants can be considered and interviewed for job position *j* during a particular amount of time in each period.

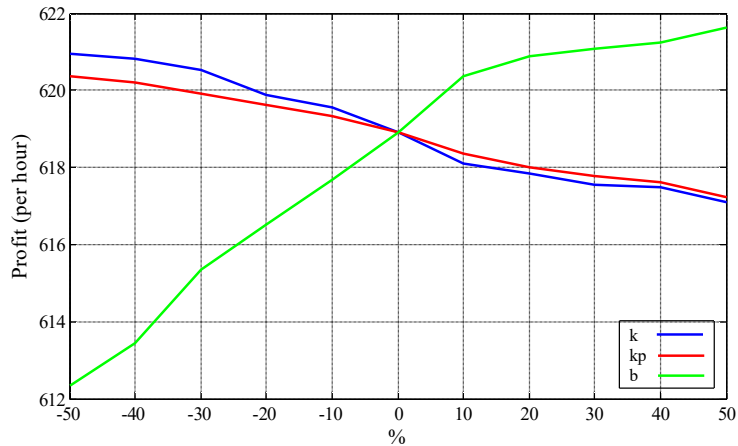

**Figure 6.** Sensitivity analysis on the mean of $\tilde{k}_j$ (blue) and $\tilde{\tilde{k}}_j$ (red) and $\tilde{b}_j$ (green).

Accordingly, to maximize the objective function, the number of hires needed for this job in period *t*, i.e., $H_j^t$, also increases (decreases). Based on Constraint (13), increasing (decreasing) $H_j^t$ causes increasing (decreasing) the corresponding growth, attrition and advancement rates in period *t*. Now, according to Constraint (14), if $Z_j^t$ and the aforementioned operational rates increase (decrease), then the number of employees for job *j* at the end of period *t*, i.e., $S_j^t$, increases (reduces). Hence, increasing (decreasing) $S_j^t$ and making $|H_j^t - Z_j^t|$ small result in increasing (reducing) the maximum total average profit generated in this department, as illustrated in Figure 6. That $H_j^t$ increase (decrease) as a result of increasing (decreasing) $Z_j^t$ implies that the cost of having plethora employees for job position *j* in period *t* does not increase after reducing the mean of either time parameter because otherwise the profit objective function does not maximize. Moreover, if Constraint (16) is binding for job position *j* in time period *t*, decreasing the mean of a corresponding time parameter might not lead to hiring more applicants for that job position in period *t*. Consequently, no increase in the number of available employees happens for job position *j* at the end of period *t*, which leads to increasing the total average profit yielded no longer. However, according to Figure 6, the later circumstance does not happen in the current case study. Furthermore, if we increased the mean of either time parameter associated with job position *j* too much (theoretically to positive infinity), then no new candidates would be hired for this job position. However, if this happened for all job positions simultaneously, i.e., we increased the mean of either time parameter of all jobs too much, then no candidates would be hired in this department and consequently the corresponding curve in Figure 6 would become horizontal. It implies that the optimal values of operational rates would be determined in such a way that no new hires would be needed, i.e., $H_j^t = 0, \forall j, t$.

Now, we study the effect of changes in the mean of acceptance rates on the value of objective function, where the acceptance rate of all job positions simultaneously varies by ±50%. If the acceptance rate of job position *j* decreases (increases), then fewer (more) candidates are hired to start that job in period *t*. In this case, following the arguments given above in case of reducing (increasing) the mean of either time parameter, the number of available employees in job *j* at the end of period *t*, i.e., $S_j^t$, becomes smaller (larger), which results in declining (increasing) the average profit yielded by employees in job *j*. Since this happens for all job positions, the total average profit generated in this department reduces (increases) as illustrated in Figure 6.

## 7. Conclusions

In the present paper, we proposed an innovative multi-period mixed integer nonlinear programming model in a stochastic environment to plan and manage the talents of an organization. The model sought to maximize the total average profit generated by employees in the whole time horizon. To mitigate much of the uncertainty associated with managerial decisions and talent management planning activities, we regarded some uncertain parameters in the developed model as stochastic following any arbitrary PDFs. A solution methodology based on linearization of nonlinear terms and chance-constrained programming was then proposed to tackle nonlinearity and uncertain parameters involved.

A review of the literature revealed that such a modeling of talent management problem is unique and could potentially open a new set of research streams to further studies in talent management planning. As identified earlier, organizations need to align talent management plans with strategic operations planning, to ensure that they have the right people with the right capabilities in place to deploy intended processes [46,47]. The proposed model was applied to the human resource needs for the logistics division of a large mining and construction company, demonstrating the realistic application of the developed mathematical model to provide useful and accurate planning directives to its SCM and HR departments. As a result of multi-periodicity of the proposed model, its output allocates required resources better during every period of planning horizon and ensures appropriate planning for successful outcomes. In the words of one manager interviewed, "We have a better sense of how many people we need to make initial contact with, interview, and recruiting to be able to fill our future pipeline of talent".

We seek to apply such models going forward to enable companies to more efficiently and effectively maximize their efforts related to HRM and TM. Viewing humans from a supply–demand perspective is an area of future research that will enable effective dynamic capabilities required to compete in complex global organizational environment [48]. In any future applications and developments of our model, the uncertain parameters, which were supposed to be stochastic in the present study, can be tackled by other uncertain programming approaches such as fuzzy programming. In fuzzy programming, the uncertain parameters are modeled via membership functions. Furthermore, a combination of both fuzzy and stochastic programming approaches in which we regard the uncertain parameters as either fuzzy-random or random-fuzzy variables can be applied to the proposed mathematical model.

**Acknowledgments:** All sources of funding of the study should be disclosed. Please clearly indicate grants that you have received in support of your research work. Clearly state if you received funds for covering the costs to publish in open access.

**Acknowledgments:** Funding for this study was provided by the Supply Chain Resource Cooperative in the Poole College of Management at North Carolina State University.

**Author Contributions:** Hadi Moheb-Alizadeh and Robert B. Handfield have contributed all sections of the paper in the same capacity.

**Conflicts of Interest:** No potential conflict of interest is reported by the authors.

## Appendix A

**Part (a):** In this part, we introduce the related functions of the Lingo 16 optimization software to model stochasticity in a programming problem, as follows:

- @SPSTGRNDV(1,RANDOM_VAR Name) for identifying the random variables, which are $\tilde{k}_j$ and $\tilde{\tilde{k}}_j$ in our stochastic program.
- @SPDIST<TYPE>(PARAM_1[,...,PARAM_N],RANDOM_VAR Name) for declaring parametric distributions such as lognormal distribution; @SPDISTLOGN(MU,SIGMA, RNDVAR), exponential distribution; @SPDISTEXPO(LAMDA,RNDVAR) and so forth.

- @SPCHANCE('Set_Name','>='|'<=',Probability) for identifying Chance-Constraint set; here Constraint (33).
- @SPSAMPSIZE (1,SIZE) for setting sample sizes in generating random numbers of random variables.

**Part (b):** In order to define the stochastic parameters $\tilde{k}_j$ and $\tilde{\tilde{k}}_j$ in the program body, the following commands are utilized:

@for(job_position(j): @SPSTGRNDV(1,k(j)));

@for(job_position(j): @SPSTGRNDV(1,kp(j)));

where k(j) and kp(j) denote $\tilde{k}_j$ and $\tilde{\tilde{k}}_j$, respectively. Furthermore, we write the following commands to determine the associated PDFs with $\tilde{k}_j$ and $\tilde{\tilde{k}}_j$; $j=1,\ldots,5$:

@SPDISTEXPO(2.6570,k(1)); @SPDISTEXPO(1.3422,k(2));

@SPDISTEXPO(1.1328,k(3)); @SPDISTEXPO(0.9961,k(4));

@SPDISTLOGN(0.777,0.521,k(5));

@SPDISTEXPO(1.2091,kp(1)); @SPDISTEXPO(0.8482,kp(2));

@SPDISTEXPO(0.7617,kp(3)); @SPDISTEXPO(0.6957,kp(4));

@SPDISTLOGN(1.019,0.467,kp(5));

In addition, the following commands characterize the stochastic Constraint (33):

@SPCHANCE('CCP_TIME','>=',0.95);

@SPCHANCE('CCP_TIME',C_TIME);

@SPSAMPSIZE(1,60);

where CCP_TIME indicates stochastic Constraint (33) in the program body. In fact, the function @SPCHANCE requires the probability of having $\sum_{j=1}^{n}\sum_{i=1}^{m}(a_{ij}^{t}\tilde{k}_j + \sum_{r_1=1}^{\bar{r}_1}\tilde{\tilde{k}}_j 2^{r_1}(\bar{x}_{ij}^{t})_{r_1}) \leq RT^t$ being greater than or equal to 0.95 ($=\alpha_2$). To accomplish this target, 60 random numbers of $\tilde{k}_j$ and $\tilde{\tilde{k}}_j$; $j=1,\ldots,5$ following their respective PDFs are generated by defining the function @SPSAMPSIZE within the program body.

*Logistics* **2017**, *1*, 5
29 of 308. Lewis, R.E.; Heckman, R.J. Talent management: A critical review. *Hum. Resour. Manag. Rev.* **2006**, *16*, 139–154.
9. Scullion, H.; Collings, D. *Global Talent Management*; Routledge: London, UK, 2011.
10. Collings, D.G.; Mellahi, K. Strategic talent management: A review and research agenda. *Hum. Resour. Manag. Rev.* **2009**, *19*, 304–313.
11. Thunnissen, M.; Boselie, P.; Fruytier, B. A review of talent management: 'infancy or adolescence?' *Int. J. Hum. Resour. Manag.* **2013**, *24*, 1744–1761.
12. Gungor, Z.; Serhadlıoglu, G.; Kesen, S.E. A fuzzy AHP approach to personnel selection problem. *Appl. Soft Comput.* **2009**, *9*, 641–646.
13. Lin, H.-T. Personnel selection using analytic network process and fuzzy data envelopment analysis approaches. *Comput. Ind. Eng.* **2010**, *59*, 937–944.
14. Gibney, R.; Shang, J. Decision making in academia: A case of the dean selection process. *Math. Comput. Model.* **2007**, *46*, 1030–1040.
15. Liao, S.-K.; Chang, K.-L. Selecting public relations personnel of hospitals by analytic network process. *J. Hosp. Mark. Public Relat.* **2009**, *19*, 52–63.
16. Kelemenis, A.; Askounis, D. A new TOPSIS-based multi-criteria approach to personnel selection. *Expert Syst. Appl.* **2010**, *37*, 4999–5008.
17. Dursun, M.; Karsak, E.E. A fuzzy MCDM approach for personnel selection. *Expert Syst. Appl.* **2010**, *37*, 4324–4330.
18. Kelemenis, A.; Ergazakis, K.; Askounis, D. Support managers' selection using an extension of fuzzy TOPSIS. *Expert Syst. Appl.* **2011**, *38*, 2774–2782.
19. Sang, X.; Liu, X.; Qin, J. An analytical solution to fuzzy TOPSIS and its application in personnel selection for knowledge-intensive enterprise. *Appl. Soft Comput.* **2015**, *30*, 190–204.
20. Kabak, M.; Burmaoglu, S.; Kazançoglu, Y. A fuzzy hybrid MCDM approach for professional selection. *Expert Syst. Appl.* **2012**, *39*, 3516–3525.
21. Hedge, J.W.; Borman, W.C.; Bourne, M.J. Designing a system for career development and advancement in the U.S. Navy. *Hum. Resour. Manag. Rev.* **2006**, *16*, 340–355.
22. Baruch, Y. Transforming careers: From linear to multidirectional career path. *Career Dev. Int.* **2004**, *9*, 58–73.
23. Hamori, M.; Koyuncu, B. Career advancement in large organizations in Europe and the United States: Do international assignments add value? *Int. J. Hum. Resour. Manag.* **2011**, *22*, 843–862.
24. Garavan, T.N.; O'Brien, F.; O'Hanlon, D. Career advancement of hotel managers since graduation: A comparative study. *Pers. Rev.* **2006**, *35*, 252–280.
25. Morrell, K.M.; Loan-Clarke, J.; Wilkinson, A.J. Organizational change and employee turnover. *Pers. Rev.* **2004**, *33*, 161–173.
26. Hom, P.W.; Kinicki, A.J. Toward a greater understanding of how dissatisfaction drives employee turnover. *Acad. Manag. J.* **2001**, *44*, 975–987.
27. Maertz, C.P.; Griffeth, R.W.; Campbell, N.S.; Allen, D.G. The effects of perceived organizational support and perceived supervisor support on employee turnover. *J. Organ. Behav.* **2007**, *28*, 1059–1075.
28. Steel, R.P.; Lounsbury, J.W. Turnover process models: Review and synthesis of a conceptual literature. *Hum. Resour. Manag. Rev.* **2009**, *19*, 271–282.
29. Holtom, B.C.; Mitchell, T.R.; Lee, T.W.; Eberly, M.B. Turnover and retention research: A glance at the past, a closer review of the present, and a venture into the future. *Acad. Manag. Ann.* **2008**, *2*, 231–274.
30. Arlotto, A.; Chick, S.E.; Gans, N. Hiring and retention policies for workers who learn. *Manag. Sci.* **2014**, *60*, 110–129.
31. Kyndt, E.; Dochy, F.; Michielsen, M.; Moeyaert, B. Employee retention: Organizational and personal perspectives. *Vocat. Learn.* **2009**, *2*, 195–215.
32. Ramlall, S. A review of employee motivation theories and their implications for employee retention within organization. *J. Am. Acad. Bus.* **2004**, *5*, 52–63.
33. Whitt, W. The impact of increased employee retention on performance in a customer contact center. *Manuf. Serv. Oper. Manag.* **2006**, *8*, 235–252.
34. Hausknecht, J.P.; Rodda, J.; Howard, M.J. Targeted employee retention: Performance-based and job-related differences in reported reasons for staying. *Hum. Resour. Manag.* **2009**, *48*, 269–288.
35. Hwang, C.L.; Yoon, K. *Multiple Attributes Decision Making Methods and Applications*; Springer: Berlin, Germany, 1981.